\documentclass[12pt]{amsart}
\usepackage{graphicx,enumerate,amsmath,amssymb,afterpage,cite}
\usepackage[sc]{mathpazo}

\newcommand\bigO{\mathcal{O}}
\newcommand{\dif}{\,\mathrm{d}}

\newcommand{\parderiv}[2]{\frac{\partial #1}{\partial #2}}
\newcommand{\deriv}[2]{\frac{\mathrm{d} #1}{\mathrm{d} #2}}
\newcommand\BS\boldsymbol
\newcommand\picturehere[1]{\includegraphics[width=0.5\textwidth]{#1}}
\newcommand\largepicturehere[1]{\includegraphics[width=\textwidth]{#1}}
\newenvironment{acknowledgment}{\section*{Acknowledgment}}{}

\numberwithin{equation}{section}
\numberwithin{figure}{section}
\numberwithin{table}{section}

\begin{document}

\title[A simple linear response closure approximation for slow dynamics]
{A simple linear response closure approximation for slow
  dynamics of a multiscale system with linear coupling}

\author{Rafail V. Abramov}

\address{Department of Mathematics, Statistics and Computer
  Science\\University of Illinois at Chicago\\851 S. Morgan st. (M/C
  249)\\ Chicago, IL 60607}

\email{abramov@math.uic.edu}


\date{\today}

\pagestyle{myheadings}

\begin{abstract}

Many applications of contemporary science involve multiscale dynamics,
which are typically characterized by the time and space scale
separation of patterns of motion, with fewer slowly evolving variables
and much larger set of faster evolving variables. This time-space
scale separation causes direct numerical simulation of the evolution
of the dynamics to be computationally expensive, due both to the large
number of variables and the necessity to choose a small discretization
time step in order to resolve the fast components of dynamics. In this
work we propose a simple method of determining the closed model for
slow variables alone, which requires only a single computation of
appropriate statistics for the fast dynamics with a certain fixed
state of the slow variables. The method is based on the first-order
Taylor expansion of the averaged coupling term with respect to the
slow variables, which can be computed using the linear
fluctuation-dissipation theorem. We show that, with simple linear
coupling in both slow and fast variables, this method produces quite
comparable statistics to what is exhibited by a complete two-scale
model. The main advantage of the method is that it applies even when
the statistics of the full multiscale model cannot be simulated due to
computational complexity, which makes it practical for real-world
large scale applications.
\end{abstract}

\maketitle

\section{Introduction}

Multiscale dynamics are common in applications of contemporary
science, such as geophysical science and climate change prediction
\cite{FraMajVan,Has,BuiMilPal,Pal3}. Multiscale dynamics are typically
characterized by the time and space scale separation of patterns of
motion, with (typically) fewer slowly evolving variables and much
larger set of faster evolving variables. This time-space scale
separation causes direct numerical simulation of the evolution of the
dynamics be computationally expensive, due both to the large number of
variables and the necessity to choose a small discretization time step
in order to resolve the fast components of dynamics.

In the climate change science the situation is further complicated by
the fact that climate is characterized by the {\em long-term
  statistics of the slow variables}, which, under small changes of
parameters (such as the solar radiation forcing, greenhouse gas
content, etc) change over even longer time scale than the motion of
the slow variables themselves. In this situation, where long-term
statistics of the slow motion patterns need to be captured, the direct
forward time integration of the most comprehensive global circulation
models (GCM) is subject to enormous computational expense.

As a more computationally feasible alternative to direct forward time
integration of the complete multiscale model, it has long been
recognized that, if a closed simplified model for the slow variables
alone is available, one could use this closed slow-variable model
instead to simulate the statistics of the slow variables. Numerous
closure schemes were developed for multiscale dynamical systems
\cite{CroVan,FatVan,MajTimVan,MajTimVan2,MajTimVan3,MajTimVan4}, which
are all based on the averaging principle over the fast variables
\cite{Pap,Van,Vol}. Some of the methods (such as those in
\cite{MajTimVan,MajTimVan2,MajTimVan3,MajTimVan4}) replace the fast
nonlinear dynamics with suitable stochastic processes \cite{Wilks} or
conditional Markov chains \cite{CroVan}, while others \cite{FatVan}
provide direct closure by suitable tabulation and curve
fitting. However, it seems that all these approaches require either
extensive computations to produce a closed model (for
example,\cite{CroVan,FatVan} require multiple simulations of fast
variables alone with different fixed states of slow variables), or
somewhat {\em ad hoc} determination of closure coefficients by
matching areas under the time correlation functions
\cite{MajTimVan,MajTimVan2,MajTimVan3,MajTimVan4}.

In this work we propose a simple method of determining the closed
model for slow variables alone, which requires only a single
computation of appropriate statistics for the fast dynamics with a
certain fixed state of the slow variables. The method is based on the
first-order Taylor expansion of the averaged coupling term with
respect to the slow variables, which, as we show, can be computed
using the linear fluctuation-dissipation theorem
\cite{Abr5,Abr6,Abr7,Abr8,AbrMaj4,AbrMaj5,AbrMaj6,AbrMaj7,MajAbrGro,Ris}.
We show through the computations with the appropriately rescaled
two-scale Lorenz 96 model \cite{Abr8} that, with simple linear
coupling in both slow and fast variables, this method produces quite
comparable statistics to what is exhibited by the slow variables of
the complete two-scale Lorenz model. The main advantage of the method
is its simplicity and easiness of implementation, partly due to the
fact that the fast dynamics need not be explicitly known (that is, the
fast dynamics can be provided as a ``black-box'' algorithm), and the
parameters of the closed model for the slow variables are determined
from the appropriate statistics of the fast variables for a given
fixed state of the slow variables. Additionally, the method can be
applied even when the statistics for slow variables of the full
multiscale model are not available due to computational expense.

The manuscript is structured as follows. In Section \ref{sec:theory}
we outline the theoretical grounds for the algorithm. In Section
\ref{sec:lorenz} we describe the two-scale Lorenz model
\cite{FatVan,Lor,LorEma}, appropriately rescaled so that the means and
variances of both the fast and slow variables are near zero and one,
respectively \cite{Abr8}. Section \ref{sec:num} contains the results
of numerical experiments with both the two-scale Lorenz model and the
reduced set of equations for slow variables only, comparing various
statistics of the time series. Section \ref{sec:sum} summarizes the
results of this work.

\section{The linear response closure approximation}
\label{sec:theory}

Consider a two-scale system of differential equations of the form
\begin{equation}
\label{eq:dyn_sys}
\deriv{\BS x}t=\BS F(\BS x,\BS y),\qquad
\deriv{\BS y}t=\BS G(\BS x,\BS y),
\end{equation}
where $\BS x=\BS x(t)\in\mathbb R^{N_x}$ are the slow variables, $\BS
y=\BS y(t)\in\mathbb R^{N_y}$ are the fast variables, and $\BS F$ and
$\BS G$ are $N_x$ and $N_y$ vector-valued functions of $\BS x$ and
$\BS y$, respectively. Here and below, we assume that the fast
variables $\BS y$ are ``unresolved'', that is, the computation of the
dynamics for $\BS y$ requires such a small time discretization step,
that the direct computation of \eqref{eq:dyn_sys} for long time
intervals is practically infeasible. We also assume that the dynamics
for $\BS y$ are fast enough for the system in \eqref{eq:dyn_sys} to be
approximated by the averaged dynamics for $\BS x$, given by
\begin{equation}
\label{eq:dyn_sys_slow_limiting_x}
\deriv{\BS x}t=\langle\BS F\rangle(\BS x),\qquad \langle\BS
F\rangle(\BS x)=\int_{\mathbb R^{N_y}}\BS F(\BS x,\BS z) \dif\mu_{\BS
  x}(\BS z),
\end{equation}
for finite times (for a more detailed description of the averaging
formalism, see \cite{Abr6,Abr8,Pap,Van,Vol}). Above, $\mu_{\BS x}$ is
the invariant probability measure of the limiting fast dynamics given
by
\begin{equation}
\label{eq:dyn_sys_fast_limiting_z}
\deriv{\BS z}\tau=\BS G(\BS x,\BS z),
\end{equation}
where the solution is given by the flow $\BS z(\tau)=\phi_{\BS
  x}^\tau\BS z_0$, while $\BS x$ is a fixed constant parameter for
\eqref{eq:dyn_sys_fast_limiting_z} and, consequently, $\mu_{\BS
  x}$. Below we assume that $\langle\BS F\rangle(\BS x)$ varies
smoothly with respect to $\BS x$ (it is known that for stochastically
driven systems this property is generic, and for deterministic
dynamics this happens when $\mu_{\BS x}$ is an SRB measure
\cite{EckRue,Rue1,Rue2,Rue3,You}). Under the ergodicity assumption for
$\mu_{\BS x}$, one can replace the measure average with time average:
\begin{equation}
\label{eq:F_time_average}
\langle\BS F\rangle(\BS x)=\lim_{r\to\infty}\frac 1r\int_0^r\BS F(\BS x,
\BS z(\tau))\dif\tau,
\end{equation}
where $\BS z(\tau)$ is a long-term trajectory of
\eqref{eq:dyn_sys_fast_limiting_z}.

While suitable multiscale methods exist (see \cite{ELi} and references
therein) where the computation of the average in
\eqref{eq:F_time_average} is performed relatively rarely, even that
might not be computationally feasible in complex multiscale systems
with many fast variables. For the purpose of this work, here we assume
that the computation of \eqref{eq:dyn_sys_fast_limiting_z} is
practically feasible only for a single choice of the constant
parameter $\BS x=\BS x^*$, where $\BS x^*$ is a suitable point, in the
vicinity of which the motion occurs, such as the mean state of the
original dynamics in \eqref{eq:dyn_sys}, or a nearby state. A poor
man's approach in this case is to compute the approximate average with
respect to $\mu_{\BS x^*}$, which is a zero order approximation:
\begin{equation}
\label{eq:F_time_average_2}
\langle\BS F\rangle(\BS x)=\int_{\mathbb R^{N_y}}\BS F(\BS x,\BS
z)\dif\mu_{\BS x^*}(\BS z)+\bigO(\|\BS x-\BS x^*\|).
\end{equation}
Here, one has to compute the time average in \eqref{eq:F_time_average}
only once, for the time series of \eqref{eq:dyn_sys_fast_limiting_z}
corresponding to $\BS x=\BS x^*$.  However, as recently found in
\cite{Abr8}, this approximation may fail to capture the chaotic
properties of the slow variables in \eqref{eq:dyn_sys}. Here we
propose the following first order correction:
\begin{equation}
\begin{split}
\langle\BS F\rangle(\BS x)&=\int_{\mathbb R^{N_y}}\BS F(\BS x,\BS
z)\dif\mu_{\BS x^*}(\BS z)+\\+&\left[\int_{\mathbb R^{N_y}}\BS F(\BS x,\BS
z)\dif\mu_{\BS x^*}^\prime(\BS z)\right](\BS x-\BS x^*)+
\bigO(\|\BS x-\BS x^*\|^2),
\end{split}
\end{equation}
where $\mu_{\BS x^*}^\prime$ is the derivative of $\mu_{\BS x}$ with
respect to $\BS x$, computed at the point $\BS x=\BS x^*$. It is shown
in \cite{Abr8} that the average of $\BS F$ with respect to $\mu_{\BS
  x^*}^\prime$ can be computed as the corresponding infinite time
linear response to infinitesimal perturbation of $\BS x$ (for details
on the linear response theory, see
\cite{Abr5,Abr6,Abr7,Abr8,AbrMaj4,AbrMaj5,AbrMaj6,AbrMaj7,MajAbrGro,Ris}
and references therein):
\begin{equation}
\label{eq:mu_deriv}
\int_{\mathbb R^{N_y}}\BS F(\BS x,\BS z)\dif\mu_{\BS x^*}^\prime(\BS
z)= \int_0^\infty\int_{\mathbb R^{N_y}} \parderiv{\BS F}{\BS y}(\BS
x,\phi_{\BS x^*}^s\BS z)\BS T_{\BS x^*,\BS z}^s\parderiv{\BS G}{\BS
  x}(\BS x^*,\BS z)\dif\mu_{\BS x^*} (\BS z)\dif s,
\end{equation}
where $\partial\BS F/\partial\BS y$ is the Jacobian with respect to
the second argument of $\BS F$, and $\BS T_{\BS x^*,\BS z}^\tau$ is
the tangent map of \eqref{eq:dyn_sys_fast_limiting_z}:
\begin{equation}
\BS T_{\BS x^*,\BS z}^\tau=\parderiv{}{\BS z}\phi_{\BS x^*}^\tau\BS z.
\end{equation}
%
Under the assumption of ergodicity of $\mu_{\BS x^*}$,
\eqref{eq:mu_deriv} can be written as the integral of the time
autocorrelation function, which is computed along a single long-time
trajectory $\BS z^*(\tau)$ of \eqref{eq:dyn_sys_fast_limiting_z} with
$\BS x=\BS x^*$:
\begin{equation}
\label{eq:mu_deriv_t}
\begin{split}
\int_{\mathbb R^{N_y}}&\BS F(\BS x,\BS z)\dif\mu_{\BS x^*}^\prime(\BS z)=\\
=&\int_0^\infty\left[\lim_{r\to\infty}\int_0^r\parderiv{\BS F}{\BS y}
(\BS x,\BS z^*(\tau+s))\BS T_{z^*(\tau)}^s\parderiv{\BS G}{\BS x}
(\BS x^*,\BS z^*(\tau))\dif\tau\right]\dif s,
\end{split}
\end{equation}
which is also done only once along the same long-term trajectory of
\eqref{eq:dyn_sys_fast_limiting_z}, as for the time average in
\eqref{eq:F_time_average_2} (for practical purposes, it can be assumed
that the time autocorrelation function decays sufficiently fast to
replace the improper integral from zero to infinity with a proper
integral up to a sufficiently long time). While formally valid, the
above formula can be unsuitable for practical computation due to the
fact that usually the limiting fast dynamics in
\eqref{eq:dyn_sys_fast_limiting_z} are strongly chaotic and mixing
with large positive Lyapunov exponents. The presence of large positive
Lyapunov exponents causes numerical instability in the computation of
the tangent map $\BS T$ of \eqref{eq:dyn_sys_fast_limiting_z} for long
response times, and, as a result, the infinite-time linear response in
\eqref{eq:mu_deriv_t} cannot be computed. Below we consider the
special setting for \eqref{eq:dyn_sys} with linear coupling (which is
common in geophysical sciences), and use the quasi-Gaussian linear
response approximation
\cite{Abr5,Abr6,Abr7,AbrMaj4,AbrMaj5,AbrMaj6,AbrMaj7,MajAbrGro} for
the practical computation of \eqref{eq:mu_deriv_t}.

\subsection{Special case with linear coupling}

Here we consider the special setting of \eqref{eq:dyn_sys} with linear
coupling between $\BS x$ and $\BS y$:
\begin{equation}
\label{eq:dyn_sys_linear}
\BS F(\BS x,\BS y)=\BS f(\BS x)+\BS L_y\BS y,\qquad
\BS G(\BS x,\BS y)=\BS g(\BS y)+\BS L_x\BS x,
\end{equation}
where $\BS f$ and $\BS g$ are nonlinear vector functions of $\BS x$
and $\BS y$, respectively, and $\BS L_x$ and $\BS L_y$ are constant
matrices of suitable sizes. For this simplified setting, observe that
$\partial\BS G/\partial\BS x=\BS L_x$, $\partial\BS F/\partial\BS
y=\BS L_y$, and, therefore, the approximate averaged system is given
by
\begin{subequations}
\begin{equation}
\label{eq:dyn_sys_averaged_linear}
\deriv{\BS x}t=\BS f(\BS x)+\BS L_y\BS{\bar z}^*+\BS L_y\BS R^*\BS
  L_x(\BS x-\BS x^*),
\end{equation}
\begin{equation}
\BS{\bar z}^*=\lim_{r\to\infty}\frac 1r\int_0^r\BS z^*(\tau)\dif\tau,
\end{equation}
\begin{equation}
\label{eq:R_tangent}
\BS R^*=\int_0^\infty\left[\lim_{r\to\infty}\frac 1r\int_0^r
\BS T_{\BS z^*(\tau)}^s\dif\tau\right]\dif s,
\end{equation}
\end{subequations}
where $\BS z(t)$ is the solution of the fast limiting system
\begin{equation}
\label{eq:dyn_sys_fast_limiting_linear_z}
\deriv{\BS z}t=\BS g(\BS z)+\BS L_x\BS x
\end{equation}
with $\BS x$ specified as a constant parameter, and $\BS z^*(t)$
corresponds to $\BS x=\BS x^*$. The formula in \eqref{eq:R_tangent} is
usually unsuitable for direct numerical computation due to rapidly
growing Lyapunov exponents at fast scales. Instead, one can use the
quasi-Gaussian linear response approximation, where
\eqref{eq:R_tangent} becomes the integral of the time autocorrelation
function under the assumption of the Gaussian invariant measure for
\eqref{eq:dyn_sys_fast_limiting_linear_z}. With linear coupling, the
measure-averaged linear response formula \eqref{eq:mu_deriv} becomes
\begin{equation}
\label{eq:mu_deriv_G}
\int_{\mathbb R^{N_y}}\BS F(\BS x,\BS z)\dif\mu_{\BS x^*}^\prime(\BS
z)= \BS L_y\left(\int_0^\infty\int_{\mathbb R^{N_y}}\BS T_{\BS x^*,\BS
  z}^s\dif\mu_{\BS x^*} (\BS z)\dif s\right)\BS L_x.
\end{equation}
Now, the Gaussian probability density is given by
\begin{equation}
\label{eq:mu_G}
p_G^*(\BS z)=(2\pi)^{-\frac{N_y}2}\det\BS\Sigma^{*-\frac
  12}\exp\left(-\frac 12(\BS z-\BS{\bar z}^*)^T \BS\Sigma^{*-1}(\BS
z-\BS{\bar z}^*)\right),
\end{equation}
where $\BS\Sigma^*$ is the long-time covariance matrix of
\eqref{eq:dyn_sys_fast_limiting_linear_z}, computed at the point $\BS
x^*$:
\begin{equation}
\BS\Sigma^*=\lim_{r\to\infty}\frac 1r\int_0^r(\BS z^*(\tau)-\BS{\bar
  z}^*)(\BS z^*(\tau)-\BS{\bar z}^*)^T\dif\tau.
\end{equation}
Recalling that $\BS T_{\BS x^*,\BS z}^s=\partial\phi_{\BS x^*}^s\BS
z/\partial\BS z$ and replacing $\dif\mu_{\BS x^*} (\BS z)=p_G^*(\BS
z)\dif\BS z$, we obtain, after integration by parts over $\BS z$,
\begin{equation}
\int_{\mathbb R^{N_y}}\BS F(\BS x,\BS z)\dif\mu_{\BS x^*}^\prime(\BS
z)= -\BS L_y\left(\int_0^\infty\int_{\mathbb R^{N_y}}\phi_{\BS x^*}^s\BS
  z\parderiv{p_G^*(\BS z)}{\BS z}\dif\BS z\dif s\right)\BS L_x.
\end{equation}
Computing the derivative of \eqref{eq:mu_G} with respect to $\BS z$,
we obtain
\begin{equation}
\int_{\mathbb R^{N_y}}\BS F(\BS x,\BS z)\dif\mu_{\BS x^*}^\prime(\BS
z)= \BS L_y\left(\int_0^\infty\int_{\mathbb R^{N_y}}\phi_{\BS x^*}^s\BS
  z(\BS z-\BS{\bar z}^*)^Tp_G^*(\BS z)\dif\BS z\dif s\right)
\BS\Sigma^{*-1}\BS L_x.
\end{equation}
Switching back to time-averaging, we obtain
\begin{equation}
\begin{split}
\int_{\mathbb R^{N_y}}\BS F(\BS x,\BS z)&\dif\mu_{\BS x^*}^\prime(\BS
z)=\\=\BS L_y&\int_0^\infty\left[\lim_{r\to\infty}\frac 1r\int_0^r
\BS z^*(\tau+s)(\BS z^*(\tau)-\BS{\bar z}^*)^T\dif\tau
\right]\dif s\;\BS\Sigma^{*-1}\BS L_x,
\end{split}
\end{equation}
and, therefore, $\BS R^*$ in \eqref{eq:dyn_sys_averaged_linear} can
now be computed as
\begin{equation}
\label{eq:R_Gaussian}
\BS R^*=\int_0^\infty\left[\lim_{r\to\infty}\frac 1r\int_0^r
\BS z^*(\tau+s)(\BS z^*(\tau)-\BS{\bar z}^*)^T\dif\tau
\right]\dif s\;\BS\Sigma^{*-1}.
\end{equation}
For details, see
\cite{Abr5,Abr6,Abr7,AbrMaj4,AbrMaj5,AbrMaj6,AbrMaj7,MajAbrGro}).

For even better precision of the linear response computation, one can
also use the blended linear response approximation
\cite{AbrMaj4,AbrMaj5,AbrMaj6,AbrMaj7}, however, in this work we do
not implement it, as it is shown that for the model and regimes
considered, the quasi-Gaussian approximation is already quite precise.

Even with the linear coupling, the function $\BS{\bar z}(\BS x)$ (the
dependence of the mean state of
\eqref{eq:dyn_sys_fast_limiting_linear_z} on $\BS x$) is not generally
linear. Thus, the validity of the linear approximation in
\eqref{eq:dyn_sys_averaged_linear} depends on the influence (or lack
thereof) of the nonlinearity of the function $\BS{\bar z}(\BS
x)$. While rigorous estimates of the validity of the linear
approximation in \eqref{eq:dyn_sys_averaged_linear} can hardly be
provided in general case, here, instead, we try to justify it by
comparing the fast limiting system in
\eqref{eq:dyn_sys_fast_limiting_linear_z} to the Ornstein-Uhlenbeck
process \cite{OrnUhl}. Consider an Ornstein-Uhlenbeck process
of the form
\begin{equation}
\label{eq:OU}
\deriv{\BS z}\tau=-\BS\Gamma(\BS z-\BS m)+\BS L_x\BS x+\BS\sigma
\deriv{\BS W_\tau}\tau,
\end{equation}
where $\BS m$ is a constant $N_y$-vector, $\BS\Gamma$ is a constant
$N_y\times N_y$ positive-definite matrix, $\BS W_t$ is a
$K$-dimensional Wiener process, $\BS\sigma$ is a constant $N_y\times K$
matrix, and $\BS x$ is, as in
\eqref{eq:dyn_sys_fast_limiting_linear_z}, is a constant
parameter. Then, it is easy to see that the difference between the
statistical mean states of \eqref{eq:OU} corresponding to $\BS x$ and
$\BS x^*$ is
\begin{equation}
\label{eq:z_diff}
\BS{\bar z}_{OU}-\BS{\bar z}_{OU}^*=\BS\Gamma^{-1}\BS L_x(\BS x-\BS
x^*),
\end{equation}
which is valid for $(\BS x-\BS x^*)$ of an arbitrary norm. At the same
time, by the regression theorem \cite{Ris}, the time correlation
function of \eqref{eq:OU} with $\BS x=\BS x^*$ is given by
\begin{equation}
\langle\BS z_{OU}^*(t+s)(\BS z_{OU}^*(t)-\BS{\bar z}_{OU}^*)^T\rangle
=\exp(-s\BS\Gamma)\BS\Sigma_{OU}^*,
\end{equation}
where $\BS\Sigma_{OU}^*$ is the covariance matrix of the
Ornstein-Uhlenbeck process in \eqref{eq:OU} for $\BS x=\BS x^*$. Thus,
according to \eqref{eq:R_Gaussian}, the infinite-time linear response
operator for \eqref{eq:OU} is computed as
\begin{equation}
\label{eq:R_OU}
\BS R_{OU}^*=\int_0^\infty\exp(-s\BS\Gamma)\dif s=\BS\Gamma^{-1}.
\end{equation}
By comparing \eqref{eq:R_OU} with \eqref{eq:z_diff}, one can see that,
for the Ornstein-Uhlenbeck process, the quasi-Gaussian linear response
formula in \eqref{eq:R_Gaussian} is exact for an arbitrarily large
perturbation $(\BS x-\BS x^*)$. Hence, if the nonlinear process in
\eqref{eq:dyn_sys_fast_limiting_linear_z} behaves statistically
similarly to the Ornstein-Uhlenbeck process in \eqref{eq:OU}, the
averaged system in \eqref{eq:dyn_sys_averaged_linear} can be expected
to behave statistically similarly to the slow part of
\eqref{eq:dyn_sys_linear}. Below we numerically test the approximation
for slow variables in the special setting with linear coupling using
the two-scale Lorenz model \cite{Abr5,Abr6,Abr8,CroVan,FatVan,Lor}.

\section{The two-scale Lorenz model}
\label{sec:lorenz}

Here we choose the two-scale forced damped Lorenz model
\cite{Abr5,Abr6,Abr8,CroVan,FatVan,Lor} for the computational study of the
dynamical properties of a two-scale slow-fast process with generic
features of climate-weather systems, such as the presence of linearly
unstable waves, strong nonlinearity, forcing, dissipation, chaos and
mixing. The two-scale forced damped Lorenz model is given by
\begin{subequations}
\label{eq:lorenz_two_scale}
\begin{equation}
\dot x_i=x_{i-1}(x_{i+1}-x_{i-2})-x_i+F_x-\frac{\lambda_y}J
\sum_{j=1}^Jy_{i,j},
\end{equation}
\begin{equation}
\dot y_{i,j}=\frac 1\varepsilon\left[y_{i,j+1}
(y_{i,j-1}-y_{i,j+2})-y_{i,j}+F_y\right]+\frac{\lambda_x}\varepsilon
x_i,
\end{equation}
\end{subequations}
where $1\leq i\leq N_x$, $1\leq j\leq J$. The following notations are
adopted above:
\begin{itemize}
\item $\BS x$ is the set of the slow variables of size $N_x$. The
  following periodic boundary conditions hold for $\BS x$:
  $x_{i+N_x}=x_i$;
\item $\BS y$ is the set of the fast variables of size $N_y=N_xJ$
  where $J$ is a positive integer. The following boundary conditions
  hold for $\BS y$: $y_{i+N_x,j}=y_{i,j}$ and $y_{i,j+J}=y_{i+1,j}$;
\item $F_x$ and $F_y$ are the constant forcing parameters;
\item $\lambda_x$ and $\lambda_y$ are the coupling parameters;
\item $\varepsilon$ is the time scale separation parameter.
\end{itemize}
Originally in \cite{CroVan,FatVan,Lor} there was no constant forcing
$F$ term in the equation for $\BS x$-variables in
\eqref{eq:lorenz_two_scale}, however, in its absence the behavior of
the $\BS y$-variables is strongly dissipative \cite{Abr5,Abr6}. Here,
as in \cite{Abr6}, we add a constant forcing $F_y$ in the right-hand
side of the second equation in \eqref{eq:lorenz_two_scale} to induce
the strongly chaotic behavior of the $\BS y$-variables with large
positive Lyapunov exponents.

\subsection{Rescaled Lorenz model}

In the Lorenz model \eqref{eq:lorenz_two_scale}, $F_x$ and $F_y$
regulate the chaos and mixing of the $\BS x$ and $\BS y$ variables,
respectively \cite{Abr5,Abr6,Abr8}. However, the mean state and mean
energy are also affected by the changes in forcing, which affects the
mean and energy trends in coupling for the fixed coupling
parameters. To adjust the effect of coupling independently of forcing,
here we rescale the Lorenz model as in \cite{MajAbrGro}. Consider the
uncoupled Lorenz model
\begin{equation}
\label{eq:L96}
\deriv{}tx_i=x_{i-1}(x_{i+1}-x_{i-2})-x_i+F
\end{equation}
with the same periodic boundary conditions as above \cite{LorEma}.
Observe that the long term statistical mean state $\bar x$ and the
standard deviation $\beta$ in \eqref{eq:L96} are the same for all
$x_i$ due to the translational invariance. Now, we rescale $\BS x$ and
$t$ as
\begin{equation}
x_i=\bar x+\beta\hat x_i,\quad t=\frac\tau\beta,
\end{equation}
where the new variables $\BS{\hat x}$ have zero mean state and unit
standard deviation. In the rescaled variables, the Lorenz model
becomes
\begin{equation}
\label{eq:L96_rescaled}
\deriv{}\tau\hat x_i=\hat x_{i-1}(\hat x_{i+1}-\hat x_{i-2})+\frac
1\beta \left[\bar x(\hat x_{i+1}-\hat x_{i-2})-\hat x_i\right]+
\frac{F-\bar x}{\beta^2},
\end{equation}
where $\bar x$ and $\beta$ are, of course, the functions of $F$. In
addition to setting the mean state and variance of $\hat x_i$ to zero
and one, respectively, due to the time rescaling the autocorrelation
functions of $\BS z$ acquire roughly identical time scaling for any
$F$ (for details, see \cite{MajAbrGro}). Here, we similarly rescale
the two-scale Lorenz model from \eqref{eq:lorenz_two_scale}:
\begin{subequations}
\label{eq:lorenz_rescaled}
\begin{equation}
\dot x_i=x_{i-1}(x_{i+1}-x_{i-2})+\frac 1{\beta_x}
\left(\bar x(x_{i+1}-x_{i-2})-x_i\right)+\frac{F_x-\bar x}
{\beta_x^2}-\frac{\lambda_y}J\sum_{j=1}^Jy_{i,j},
\end{equation}
\begin{equation}
\begin{split}
\dot y_{i,j}=\frac 1\varepsilon\bigg[y_{i,j+1}
(y_{i,j-1}-y_{i,j+2})+\frac 1{\beta_y}
\left(\bar y(y_{i,j-1}-y_{i,j+2})-y_{i,j}\right)+\\+\frac{F_y-\bar y}
{\beta_y^2}\bigg]+\frac{\lambda_x}\varepsilon x_i,
\end{split}
\end{equation}
\end{subequations}
where $\bar x$, $\bar y$, $\beta_x$ and $\beta_y$ are the long term
means and standard deviations of the corresponding uncoupled system in
\eqref{eq:L96} with either $F_x$ or $F_y$ set as a constant
forcing. In the rescaled Lorenz model \eqref{eq:lorenz_rescaled}, the
values of $F_x$ and $F_y$ do not significantly affect the mean state,
mean energy and the time scale for both the slow variables $\BS x$ and
fast variables $\BS y$, and mostly regulate the mixing and
non-Gaussianity of the probability distributions of the long-term time
series.

The Lorenz model \eqref{eq:lorenz_rescaled} represents the setting
with linear coupling as in \eqref{eq:dyn_sys_linear}, with $\BS f$,
$\BS g$, $\BS L_x$ and $\BS L_y$ given by
\begin{subequations}
\begin{equation}
f_i(\BS x)=x_{i-1}(x_{i+1}-x_{i-2})+\frac 1{\beta_x}
\left(\bar x(x_{i+1}-x_{i-2})-x_i\right)+\frac{F_x-\bar x}
{\beta_x^2},
\end{equation}
\begin{equation}
\begin{split}
g_{i,j}(\BS y)=\frac 1\varepsilon\bigg[y_{i,j+1}
(y_{i,j-1}-y_{i,j+2})+\frac 1{\beta_y}
\left(\bar y(y_{i,j-1}-y_{i,j+2})-y_{i,j}\right)+\\+\frac{F_y-\bar y}
{\beta_y^2}\bigg],
\end{split}
\end{equation}
\begin{equation}
 \BS L_x=\frac{\lambda_x}\varepsilon\BS L^T,\qquad\BS L_y=
-\frac{\lambda_y}J\BS L,\qquad\BS L_{i,(j,k)}=\delta_{ij}.
\end{equation}
\end{subequations}
Now, one can immediately see that for the Lorenz model in
\eqref{eq:lorenz_rescaled}, the limiting system in
\eqref{eq:dyn_sys_fast_limiting_z} is given by
\begin{equation}
\deriv{\BS z}\tau=\BS g(\BS z)+\frac{\lambda_x}\varepsilon\BS L^T\BS x,
\end{equation}
while the approximate averaged system around the point $\BS x^*$ is
given by
\begin{equation}
\deriv{\BS x}t=\BS f(\BS x)-\frac{\lambda_y}J\BS L\BS{\bar z}^*
-\frac{\lambda_x\lambda_y}{\varepsilon J}\BS L\BS R^*\BS L^T
(\BS x-\BS x^*).
\end{equation}
Above, the parameter $\varepsilon$ in the denominator suggests that
the first-order correction term is of order $\varepsilon^{-1}$, which
is misleading, because $\BS R^*$ scales proportionally to
$\varepsilon$ (since it is the integral of the time autocorrelation
function with rate of decay proportional to $\varepsilon^{-1}$). For
simplicity of computation, for the Lorenz model we rescale the fast
limiting dynamics by $\varepsilon$ to bring the time scale to order 1
since the scaling parameter $\varepsilon$ is known
explicitly\footnote{Generally, when no scaling parameter is available
  explicitly in \eqref{eq:dyn_sys}, one can still multiply the
  right-hand side of \eqref{eq:dyn_sys_fast_limiting_z} by a heuristic
  small parameter to bring the time scale of
  \eqref{eq:dyn_sys_fast_limiting_z} to order 1, since the averaged
  dynamics are invariant with respect to the $\varepsilon$-scaling of
  the limiting fast dynamics (for details, see \cite{Abr8}). Or,
  equivalently, use appropriately small time step and averaging window
  for \eqref{eq:dyn_sys_fast_limiting_z}.}, which yields
\begin{equation}
\label{eq:dyn_sys_fast_limiting_z_lor}
\deriv{\BS z}\tau=\varepsilon\BS g(\BS z)+\lambda_x\BS L^T\BS x,
\end{equation}
for the fast limiting dynamics for fixed $\BS x$, while for the
averaged dynamics we obtain
\begin{equation}
\label{eq:dyn_sys_slow_limiting_x_lor}
\deriv{\BS x}t=\BS f(\BS x)-\frac{\lambda_y}J\BS L\BS{\bar z}^*
-\frac{\lambda_x\lambda_y}J\BS L\BS R^*\BS L^T(\BS x-\BS x^*),
\end{equation}
where $\BS R^*$ is computed from the time series of
\eqref{eq:dyn_sys_fast_limiting_z_lor} using the quasi-Gaussian
approximation in \eqref{eq:R_Gaussian}.

\section{Numerical experiments}
\label{sec:num}

Here we present a numerical study of the proposed approximation for
slow dynamics, applied to the rescaled Lorenz model in
\eqref{eq:lorenz_rescaled}. We compare the statistical properties of
the slow variables for the three following systems:
\begin{enumerate}
\item The complete rescaled Lorenz system from
  \eqref{eq:lorenz_rescaled};
\item The approximation for slow dynamics alone from
  \eqref{eq:dyn_sys_slow_limiting_x_lor};
\item The poor man's version of \eqref{eq:dyn_sys_slow_limiting_x_lor}
  with the first-order correction term $\BS R^*$ set to zero (further
  referred to as the ``zero-order'' system).
\end{enumerate}
The fixed parameter $\BS x^*$ for the computation of $\BS R^*$ was set
to the long-term mean state $\BS{\bar x}$ of
\eqref{eq:lorenz_rescaled} (in practical situations, a rough estimate
could be used). The quasi-Gaussian approximation in
\eqref{eq:R_Gaussian} is used to compute the first order correction
term $\BS R^*$. While it is practically impossible to compute the
improper integrals in \eqref{eq:R_Gaussian} for infinite upper limit,
in practice we use sufficiently long (but finite) limits of
integration. In particular, for all computational results presented
below, the correction term $\BS R^*$ is computed numerically as
\begin{equation}
\label{eq:R_Gaussian_2}
\BS R^*=\frac 1{T_{av}}\int_0^{T_{corr}}\int_0^{T_{av}}\BS z^*(\tau+s)
(\BS z^*(\tau)-\BS{\bar z}^*)^T\dif\tau\dif s\;\BS\Sigma^{*-1},
\end{equation}
where the averaging time window $T_{av}$ equals 10000 time units,
while the correlation time window $T_{corr}$ equals 50 time units (it
was observed that the time autocorrelation function in
\eqref{eq:R_Gaussian_2} decays essentially to zero within the 50
time-unit window for all studied regimes). The mean state $\BS{\bar
  z}^*$ and the covariance matrix $\BS\Sigma^*$ are also computed by
time-averaging with the same averaging window of 10000 time units.

\begin{figure}
\picturehere{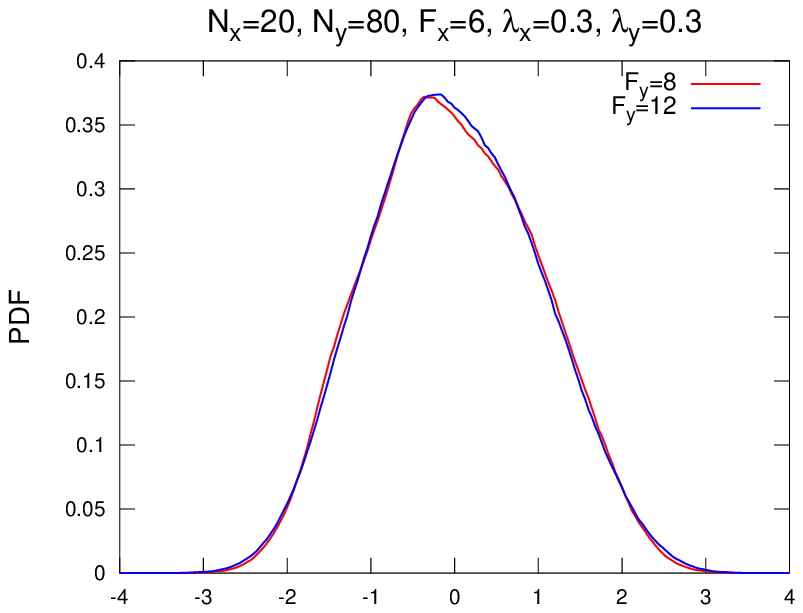}%
\picturehere{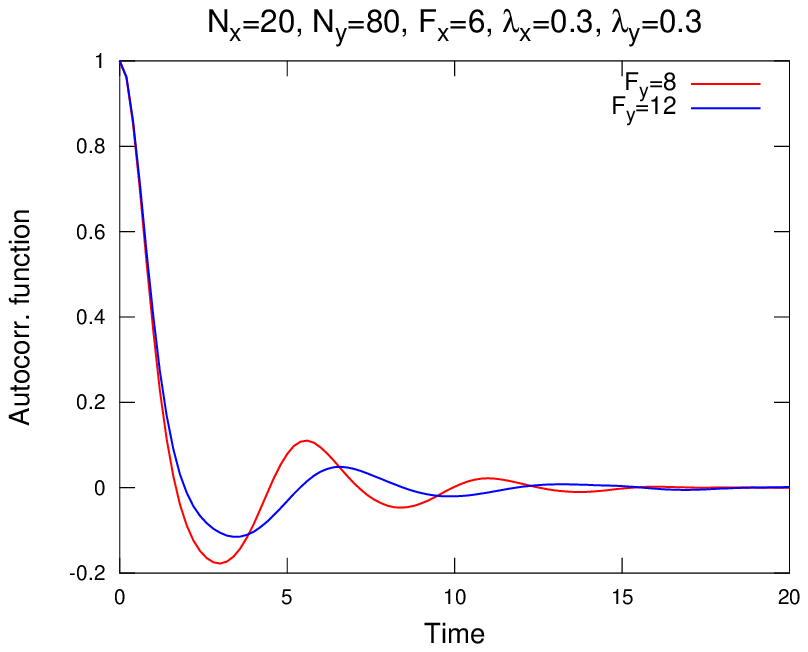}%
\caption{The probability density functions and time autocorrelation
  functions of the limiting fast dynamics in
  \eqref{eq:dyn_sys_fast_limiting_z_lor} for $F_x=6$, and
  $\lambda_x=\lambda_y=0.3$.}
\label{fig:fast_pdf_corr}
\end{figure}

Due to translational invariance of the studied models, the statistics
are invariant with respect to the index shift for the variables
$x_i$. For diagnostics, we monitor the following long-term statistical
quantities of $x_i$:
\begin{enumerate}[a.]
\item The probability density functions (PDF), computed by
  bin-counting. A PDF gives the most complete information about the
  one-point statistics of $x_i$, as it shows the statistical
  distribution of $x_i$ in the phase space.
\item The time autocorrelation functions $\langle
  x_i(t)x_i(t+s)\rangle$, where the time average is over $t$,
  normalized by the variance $\langle x_i^2\rangle$ (so that it always
  starts with 1).
\item The time cross-correlation functions $\langle
  x_i(t)x_{i+1}(t+s)\rangle$, also normalized by the variance $\langle
  x_i^2\rangle$.
\item The energy autocorrelation function $$K(s)=\frac{\langle
  x_i^2(t)x_i^2(t+s)\rangle}{\langle x_i^2\rangle^2+2\langle
  x_i(t)x_i(t+s)\rangle^2}.$$ This energy autocorrelation function
  measures the non-Gaussianity of the process (it is identically 1 for
  all $s$ if the process is Gaussian, such as the Ornstein-Uhlenbeck
  process). For details, see \cite{MajTimVan4}.
\end{enumerate}
The success (or failure) of the proposed approximation of the slow
dynamics depends on several factors. First, as the quasi-Gaussian
linear response formula is used for the computation of $\BS R^*$, the
precision will be affected by the non-Gaussianity of the fast
dynamics. Second, it depends how linearly the mean state $\BS{\bar z}$
for the fast variables depends on the slow variables $\BS x$. Here we
observe the limitations of the proposed approximation by studying a
variety of dynamical regimes of the rescaled Lorenz model in
\eqref{eq:lorenz_rescaled}. The following dynamical regimes are
studied:
\begin{itemize}
\item $N_x=20$, $J=4$ (so that $N_y=80$). Thus, the number of the fast
  variables is four times greater than the number of the slow variables.
\item $\varepsilon=0.01$. The time scale separation of two orders of
  magnitude is consistent with typical real-world geophysical
  processes (for example, the annual and diurnal cycles of the Earth's
  atmosphere).
\item $\lambda_x=\lambda_y=0.3,0.4$. These values of coupling are
  chosen so that they are neither too weak, nor too strong (although
  0.3 is weaker, and 0.4 is stronger). Recall that the standard
  deviations of both $x_i$ and $y_{i,j}$ variables are approximately
  1, and, thus, the contribution to the right-hand side from coupled
  variables is weaker than the self-contribution, but still of the
  same order.
\item $F_x=6,16$. The slow forcing $F_x$ adjusts the chaos and mixing
  properties of the slow variables, and in this work it is set to a
  weakly chaotic regime $F_x=6$, and strongly chaotic regime $F_x=16$.
\item $F_y=8,12$. The fast forcing adjusts the chaos and mixing
  properties of the fast variables. Here the value of $F_y$ is chosen
  so that the fast variables are either moderately chaotic for
  $F_y=8$, or more strongly chaotic for $F_y=12$.
\end{itemize}
In Figure \ref{fig:fast_pdf_corr} we show the probability density
functions and time autocorrelation functions for the limiting fast
dynamics in \eqref{eq:dyn_sys_fast_limiting_z_lor}, with the
parameters $N_x=20$, $N_y=80$, $F_x=6$, $\lambda_x=\lambda_y=0.3$, and
two values of the fast forcing: $F_y=8$ and $F_y=12$. Observe that the
PDFs are not Gaussian (although close to it), and have nonzero
skewness. The time autocorrelation functions decay slower for $F_y=8$
and faster for $F_y=12$, indicating slower and faster mixing,
respectively. In other regimes, these PDFs and autocorrelation
functions look very similar to what is presented in Figure
\ref{fig:fast_pdf_corr}.

\begin{figure}
\largepicturehere{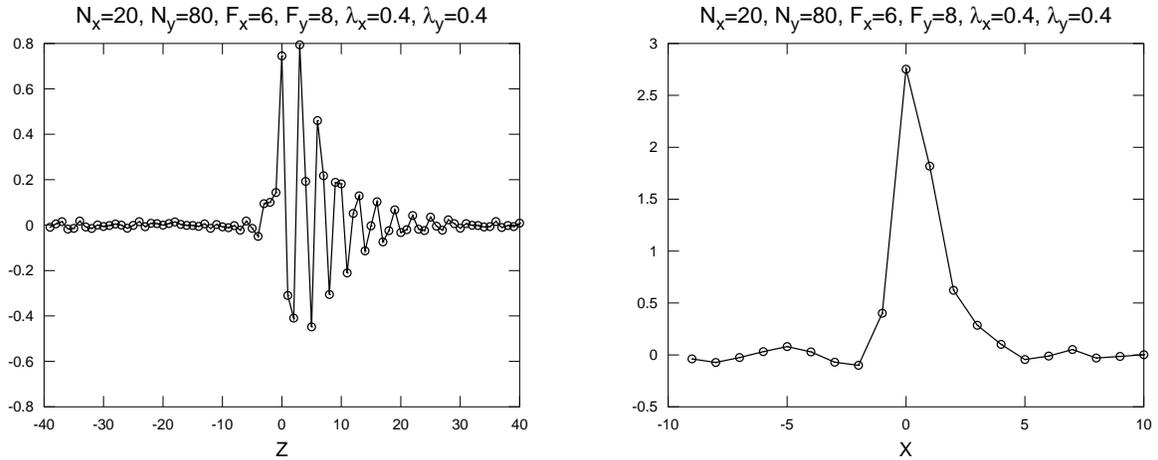}%
\caption{The matrices $\BS R^*$ (left) and $\BS L\BS R^*\BS L^T$
  (right) of the reduced dynamics in
  \eqref{eq:dyn_sys_slow_limiting_x_lor} for $F_x=6$, $F_y=8$, and
  $\lambda_x=\lambda_y=0.4$. Only a single column of each matrix is
  displayed with its diagonal element corresponding to zero horizontal
  coordinate, as the matrices are translation-invariant.}
\label{fig:resp}
\end{figure}

Another point we would like to emphasize before presenting the results
of the computational study, is that the matrices $\BS R^*$ and,
consequently, $\BS L\BS R\BS L^T$ in
\eqref{eq:dyn_sys_slow_limiting_x_lor} are not diagonal. In Figure
\ref{fig:resp} we display both $\BS R^*$ and $\BS L\BS R\BS L^T$ for
the dynamical regime with $F_x=6$, $F_y=8$, and
$\lambda_x=\lambda_y=0.4$ (only the central columns of $\BS R^*$ and
$\BS L\BS R\BS L^T$ are displayed with diagonal elements corresponding
to zero horizontal coordinates of the plots, as both matrices are
translation-invariant). Observe that there are significant
off-diagonal entries in both matrices. Both $\BS R^*$ and $\BS L\BS
R\BS L^T$ are positive-definite in the presented regime (the lowest
eigenvalues of their symmetric parts are $6.314\cdot 10^{-2}$ and
$0.8814$, respectively), and, as a result, the linear correction term
in \eqref{eq:dyn_sys_slow_limiting_x_lor} causes damping effect on the
reduced dynamics (for more details about chaotic properties of reduced
dynamics, see \cite{Abr8}). For other regimes, the matrices are
similar to the presented regime (that is, substantial off-diagonal
entries are present), and we do not display those here.

\subsection{Probability density functions of the slow dynamics}

\begin{figure}
\picturehere{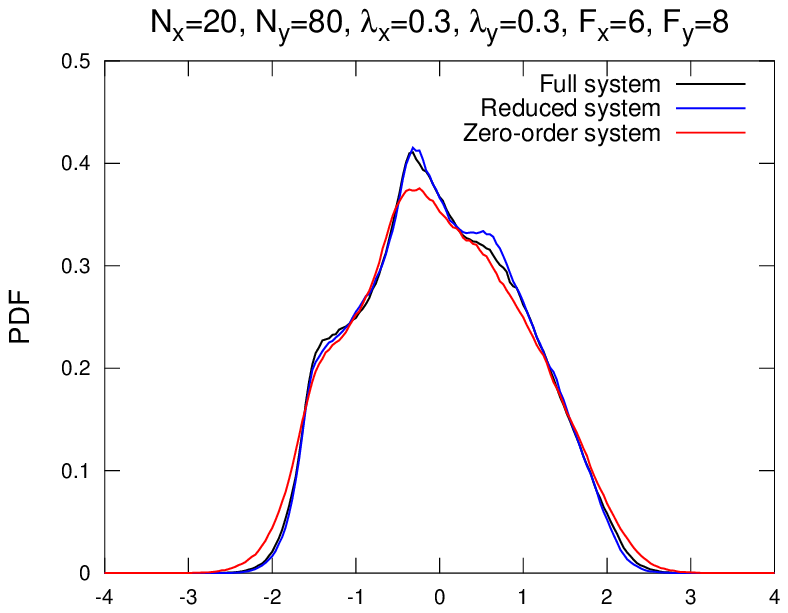}%
\picturehere{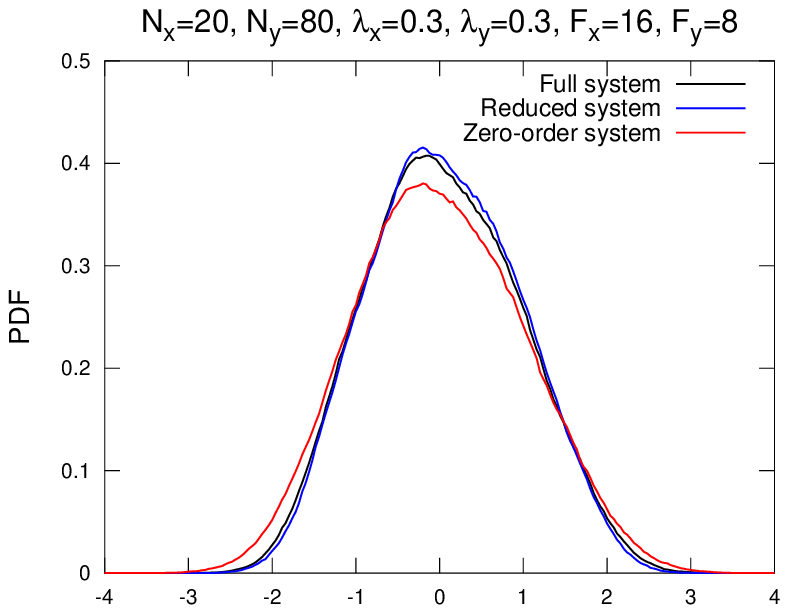}\\%
\picturehere{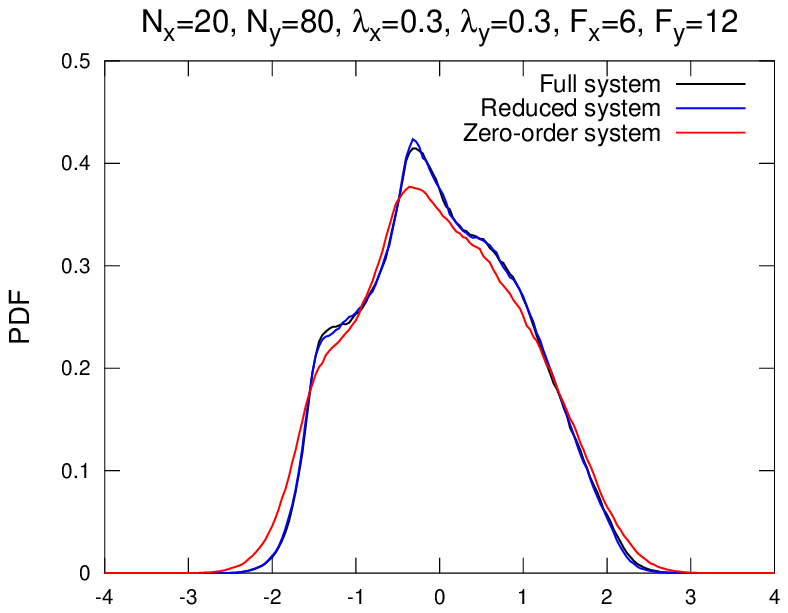}%
\picturehere{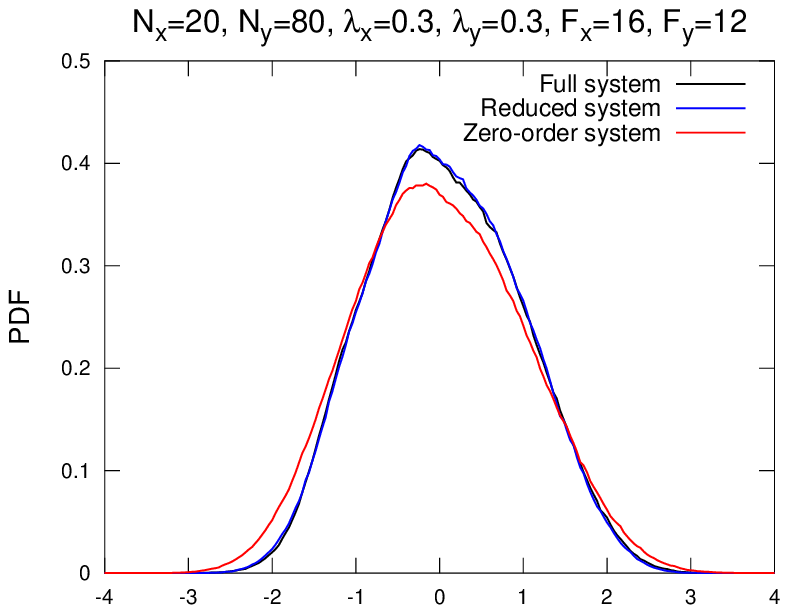}%
\caption{PDFs, $\lambda_x=\lambda_y=0.3$.}
\label{fig:pdf_l03}
\end{figure}

\begin{figure}
\picturehere{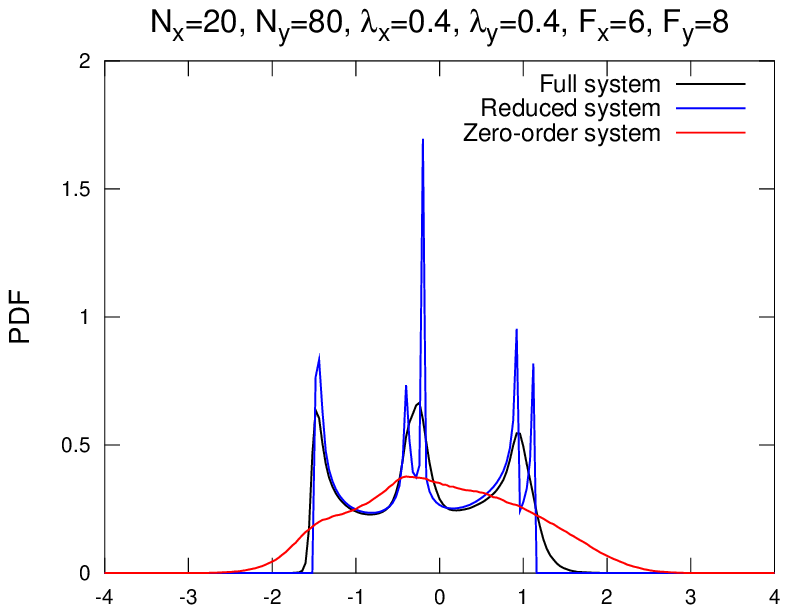}%
\picturehere{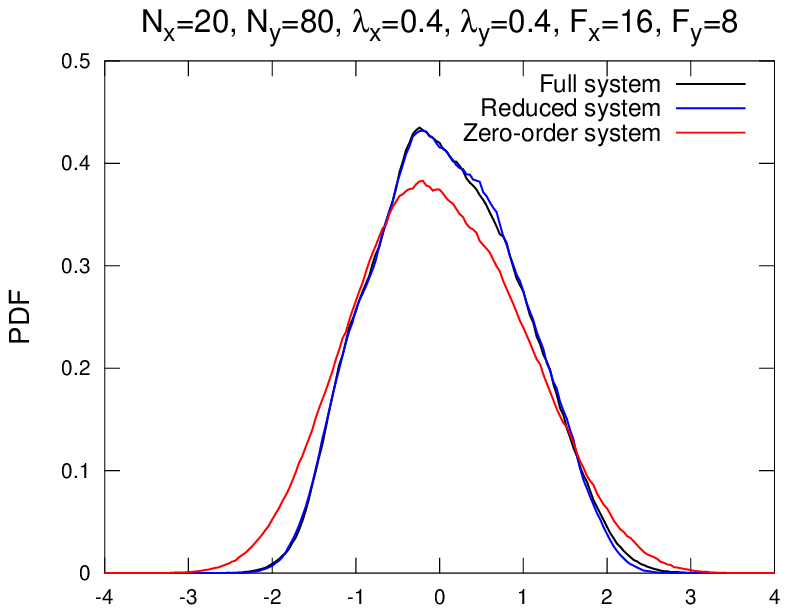}\\%
\picturehere{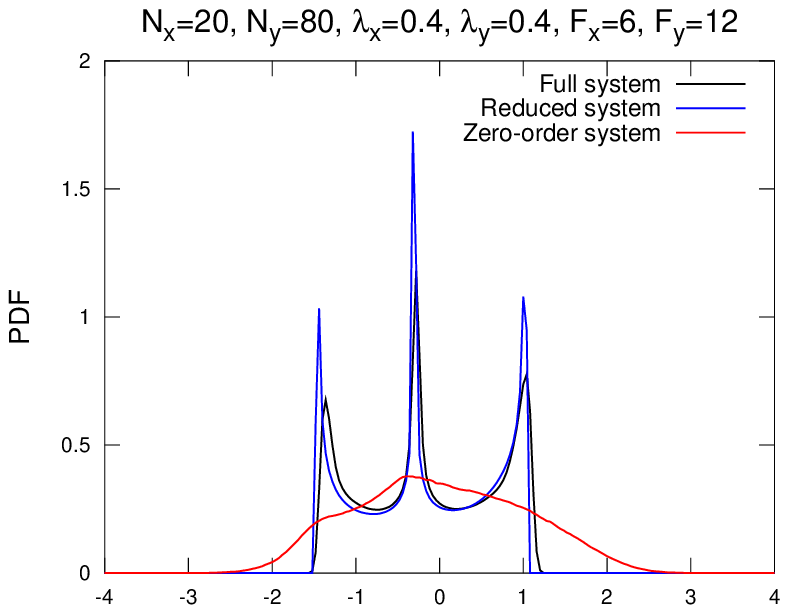}%
\picturehere{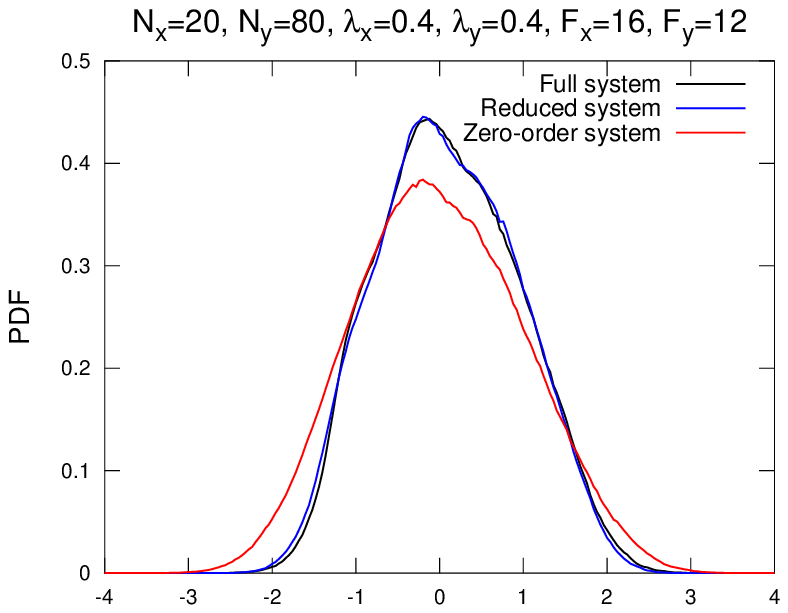}%
\caption{PDFs, $\lambda_x=\lambda_y=0.4$.}
\label{fig:pdf_l04}
\end{figure}

In Figures \ref{fig:pdf_l03} and \ref{fig:pdf_l04} we show the
probability density functions of the slow dynamics for the full
two-scale Lorenz model, the reduced closed model for the slow
variables alone in \eqref{eq:dyn_sys_slow_limiting_x_lor}, and its
poor man's zero order version without the linear correction
term. Observe that for the more weakly coupled regimes with
$\lambda_x=\lambda_y=0.3$ the PDFs look rather similar, however, it
can be seen that the reduced model with the correction term reproduces
the PDFs much closer to those of the full two-scale Lorenz model, than
the zero-order model. In the more strongly coupled regime with
$\lambda_x=\lambda_y=0.4$ the situation tilts even more in favor of
the reduced model with linear correction term in
\eqref{eq:dyn_sys_slow_limiting_x_lor}: observe that for the weakly
chaotic regime with $F_x=6$ the PDFs of the full two-scale Lorenz
model have three sharp peaks, indicating strong non-Gaussianity. The
reduced model in \eqref{eq:dyn_sys_slow_limiting_x_lor} reproduces
these peaks, while its zero-order version fails.
In addition, in Table \ref{tab:pdfs} we
show the $L_2$-errors in PDFs between the full two-scale Lorenz model
and the two reduced models. Observe that, generally, the reduced
system with the linear correction term in
\eqref{eq:dyn_sys_slow_limiting_x_lor} produces more precise results
than its poor man's version without the correction term.
\begin{table}
\begin{center}
\begin{tabular}{|c|}
\hline
\begin{tabular}{c|c}
$\lambda_{x,y}=0.3$, $F_y=8$ & $\lambda_{x,y}=0.3$, $F_y=12$ \\
\hline
\begin{tabular}{c||c|c}
 & Red. & Z.O. \\
\hline
$F_x=6$ & $5.036\cdot 10^{-3}$ & $1.165\cdot 10^{-2}$ \\
$F_x=16$ & $5.593\cdot 10^{-3}$ & $1.469\cdot 10^{-2}$ \\
\end{tabular}
&
\begin{tabular}{c||c|c}
 & Red. & Z.O. \\
\hline
$F_x=6$ & $2.581\cdot 10^{-3}$ & $1.576\cdot 10^{-2}$ \\
$F_x=16$ & $2.71\cdot 10^{-3}$ & $1.818\cdot 10^{-2}$ \\
\end{tabular}\\
\end{tabular}\\
\hline
\begin{tabular}{c|c}
$\lambda_{x,y}=0.4$, $F_y=8$ & $\lambda_{x,y}=0.4$, $F_y=12$ \\
\hline
\begin{tabular}{c||c|c}
 & Red. & Z.O. \\
\hline
$F_x=6$ & $0.1022$ & $8.857\cdot 10^{-2}$ \\
$F_x=16$ & $3.725\cdot 10^{-3}$ & $2.703\cdot 10^{-2}$ \\
\end{tabular}
&
\begin{tabular}{c||c|c}
 & Red. & Z.O. \\
\hline
$F_x=6$ & $9.28\cdot 10^{-2}$ & $0.1113$ \\
$F_x=16$ & $5.885\cdot 10^{-3}$ & $3.209\cdot 10^{-2}$ \\
\end{tabular}\\
\end{tabular}\\
\hline
\end{tabular}
\end{center}
\caption{$L_2$-errors between the PDFs of the slow variables of the
  full two-scale Lorenz model and the two reduced models. Notations:
  ``Red.'' stands for ``Reduced'' (that is,
  \eqref{eq:dyn_sys_slow_limiting_x_lor}), and ``Z.O.'' stands for
  ``Zero-order'', the poor man's version of
  \eqref{eq:dyn_sys_slow_limiting_x_lor}.}
\label{tab:pdfs}
\end{table}

\subsection{Time autocorrelation functions of the slow dynamics}

\begin{figure}
\picturehere{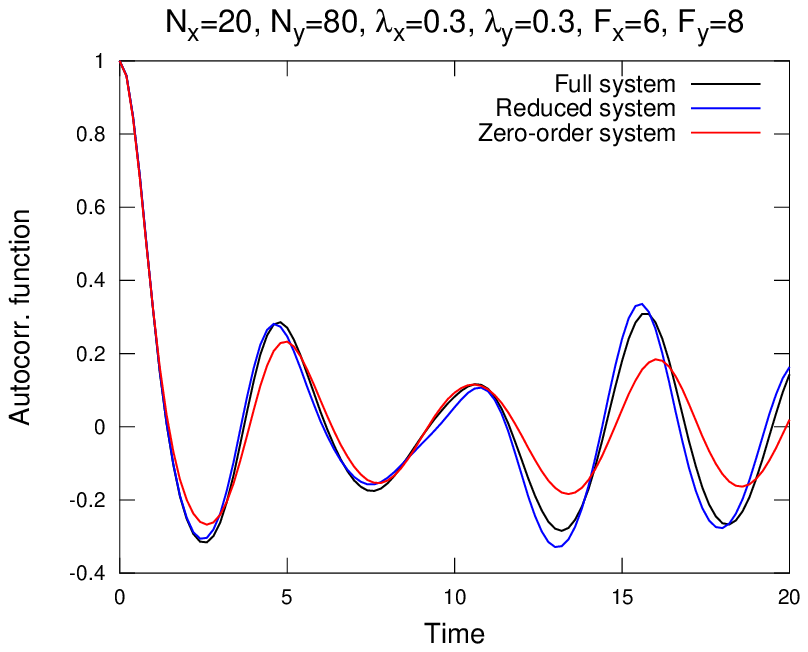}%
\picturehere{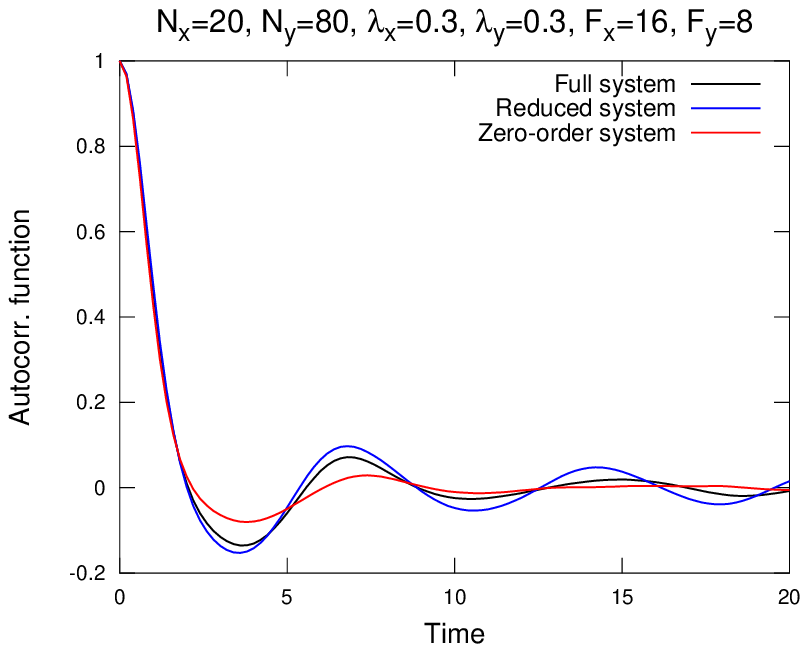}\\%
\picturehere{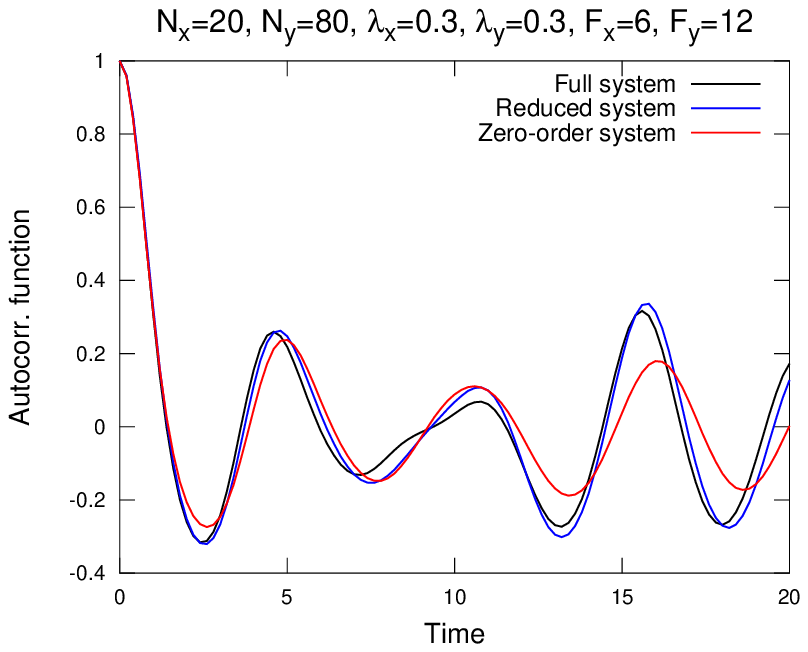}%
\picturehere{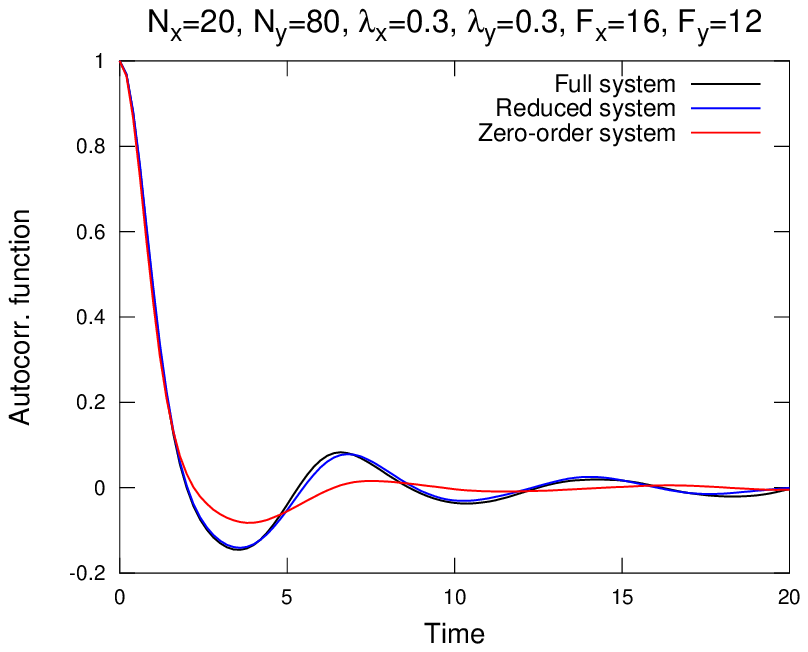}%
\caption{Time autocorrelation functions, $\lambda_x=\lambda_y=0.3$.}
\label{fig:corr_l03}
\end{figure}

\begin{figure}
\picturehere{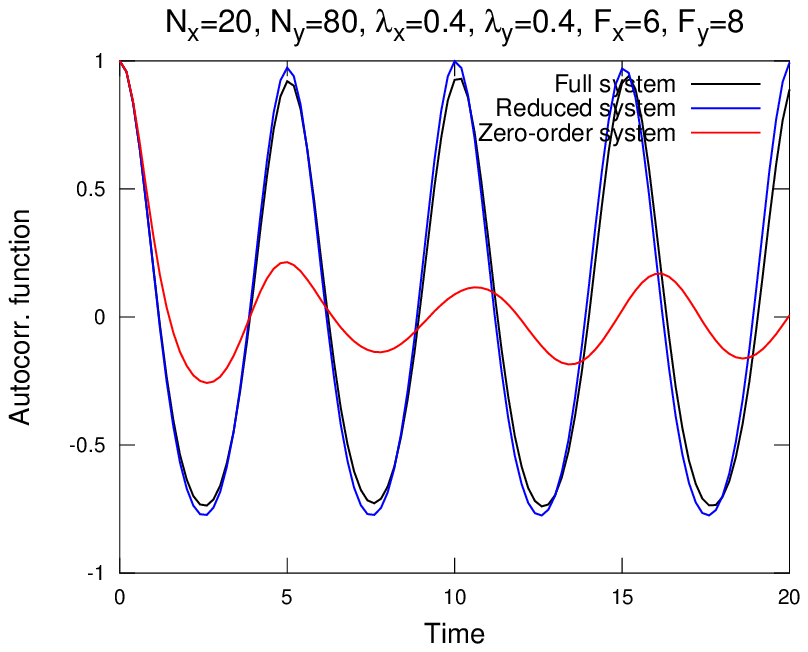}%
\picturehere{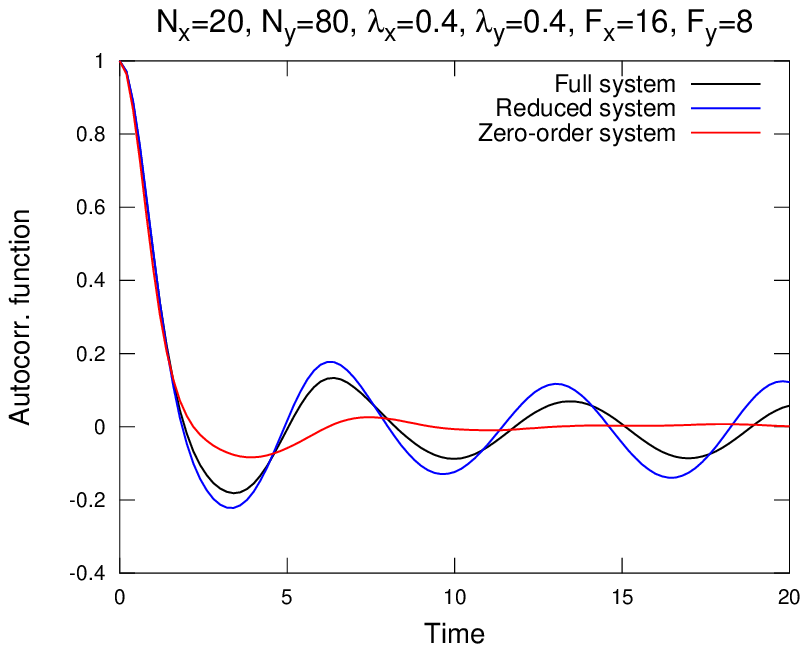}\\%
\picturehere{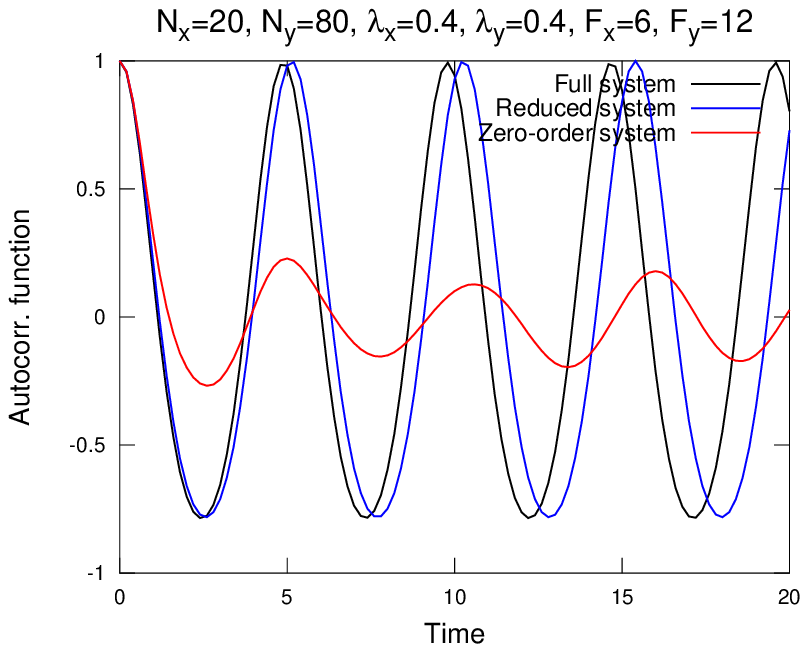}%
\picturehere{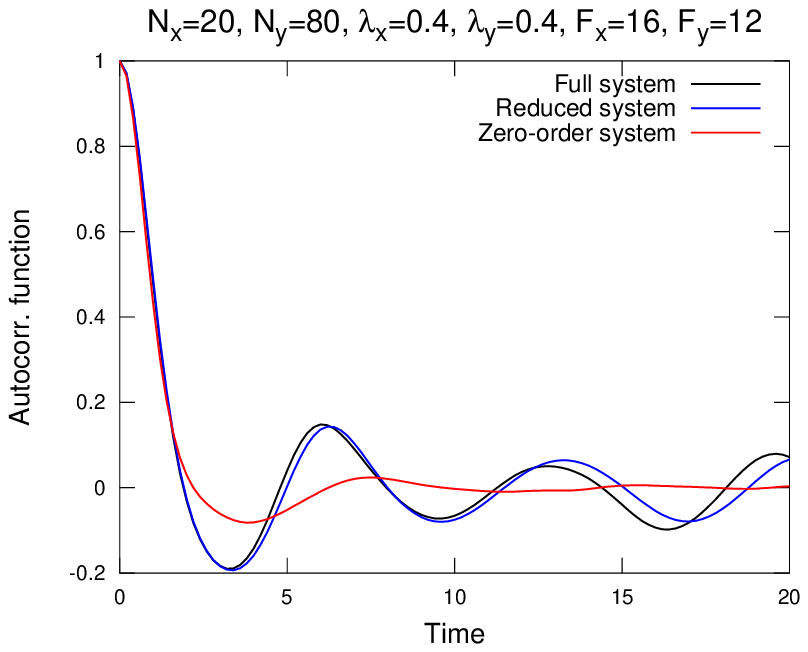}%
\caption{Time autocorrelation functions, $\lambda_x=\lambda_y=0.4$.}
\label{fig:corr_l04}
\end{figure}

In Figures \ref{fig:corr_l03} and \ref{fig:corr_l04} we show the time
autocorrelation functions of the slow dynamics for the full two-scale
Lorenz model, the reduced closed model for the slow variables alone in
\eqref{eq:dyn_sys_slow_limiting_x_lor}, and its poor man's zero order
version without the linear correction term. Observe that for the more
weakly coupled regimes with $\lambda_x=\lambda_y=0.3$ the time
autocorrelation functions look similar, yet the reduced model with the
correction term reproduces the time autocorrelation functions more
precisely than the zero-order model. In the more strongly coupled
regime with $\lambda_x=\lambda_y=0.4$ the difference between the
reduced model in \eqref{eq:dyn_sys_slow_limiting_x_lor} and its poor
man's zero-order version is even more drastic: observe that for the
weakly chaotic regime with $F_x=6$ the time autocorrelation functions
of the full two-scale Lorenz model do not exhibit decay (indicating
very weak mixing), and the reduced model in
\eqref{eq:dyn_sys_slow_limiting_x_lor} reproduces the autocorrelation
functions of the full two-scale Lorenz model rather well, while its
zero-order version fails.
In addition, in Table \ref{tab:corrs} we show the $L_2$-errors in time
autocorrelation functions (for the correlation time interval of 20
time units, as in Figures \ref{fig:corr_l03} and \ref{fig:corr_l04})
between the full two-scale Lorenz model and the two reduced
models. Observe that, generally, the reduced system with the linear
correction term in \eqref{eq:dyn_sys_slow_limiting_x_lor} produces
more precise results than its poor man's version without the
correction term.
\begin{table}
\begin{center}
\begin{tabular}{|c|}
\hline
\begin{tabular}{c|c}
$\lambda_{x,y}=0.3$, $F_y=8$ & $\lambda_{x,y}=0.3$, $F_y=12$ \\
\hline
\begin{tabular}{c||c|c}
 & Av. & Z.O. \\
\hline
$F_x=6$ & $5.841\cdot 10^{-2}$ & $0.1211$ \\
$F_x=16$ & $4.079\cdot 10^{-2}$ & $5.342\cdot 10^{-2}$ \\
\end{tabular}
&
\begin{tabular}{c||c|c}
 & Av. & Z.O. \\
\hline
$F_x=6$ & $6.539\cdot 10^{-2}$ & $0.1572$ \\
$F_x=16$ & $1.559\cdot 10^{-2}$ & $7.396\cdot 10^{-2}$ \\
\end{tabular}\\
\end{tabular}\\
\hline
\begin{tabular}{c|c}
$\lambda_{x,y}=0.4$, $F_y=8$ & $\lambda_{x,y}=0.4$, $F_y=12$ \\
\hline
\begin{tabular}{c||c|c}
 & Av. & Z.O. \\
\hline
$F_x=6$ & $5.538\cdot 10^{-2}$ & $0.3677$ \\
$F_x=16$ & $8.534\cdot 10^{-2}$ & $0.1355$ \\
\end{tabular}
&
\begin{tabular}{c||c|c}
 & Av. & Z.O. \\
\hline
$F_x=6$ & $0.2981$ & $0.3986$ \\
$F_x=16$ & $4.835\cdot 10^{-2}$ & $0.1482$ \\
\end{tabular}\\
\end{tabular}\\
\hline
\end{tabular}
\end{center}
\caption{$L_2$-errors between the time autocorrelation functions of
  the slow variables of the full two-scale Lorenz model and the two
  reduced models. Notations: ``Red.'' stands for ``Reduced'' (that is,
  \eqref{eq:dyn_sys_slow_limiting_x_lor}), and ``Z.O.'' stands for
  ``Zero-order'', the poor man's version of
  \eqref{eq:dyn_sys_slow_limiting_x_lor}.}
\label{tab:corrs}
\end{table}

\subsection{Time cross-correlation functions}

\begin{figure}
\picturehere{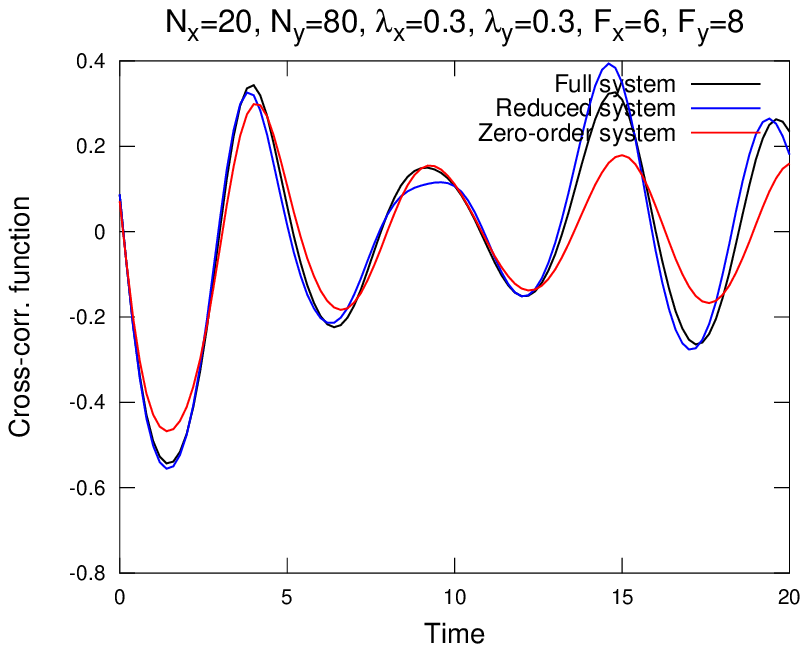}%
\picturehere{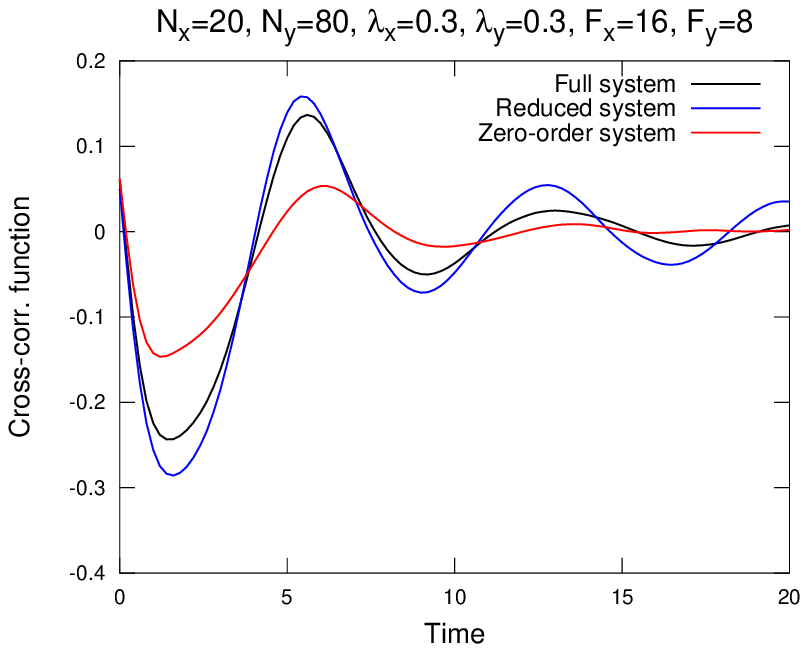}\\%
\picturehere{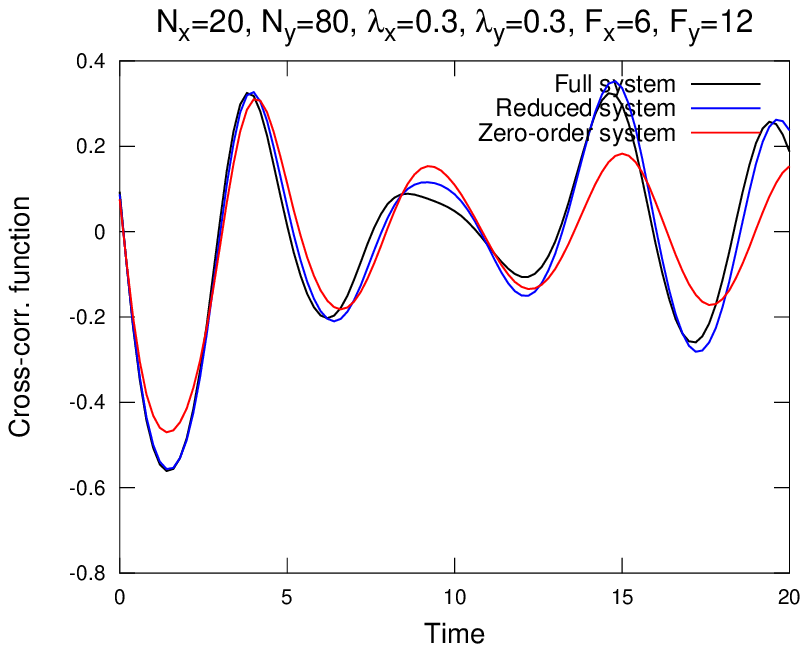}%
\picturehere{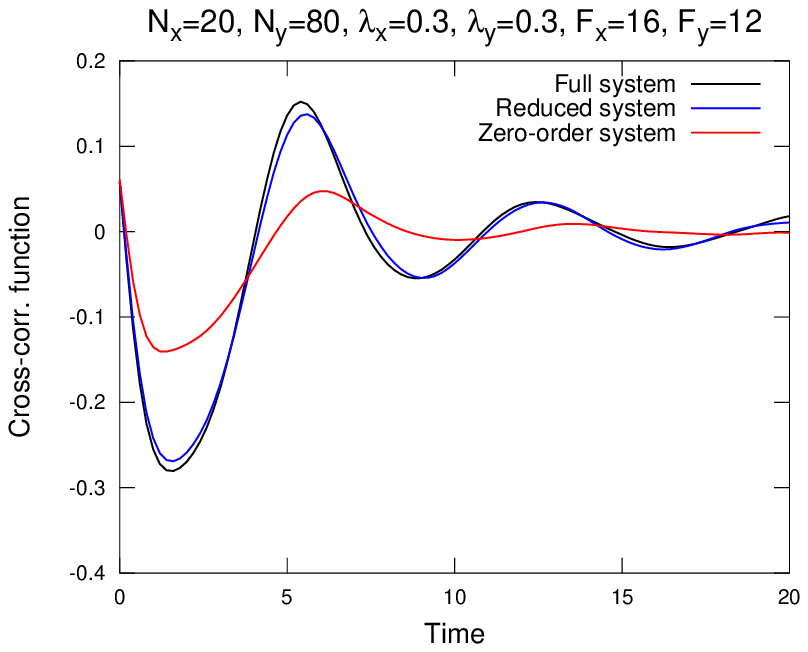}%
\caption{Time cross-correlation functions, $\lambda_x=\lambda_y=0.3$.}
\label{fig:ccorr_l03}
\end{figure}

\begin{figure}
\picturehere{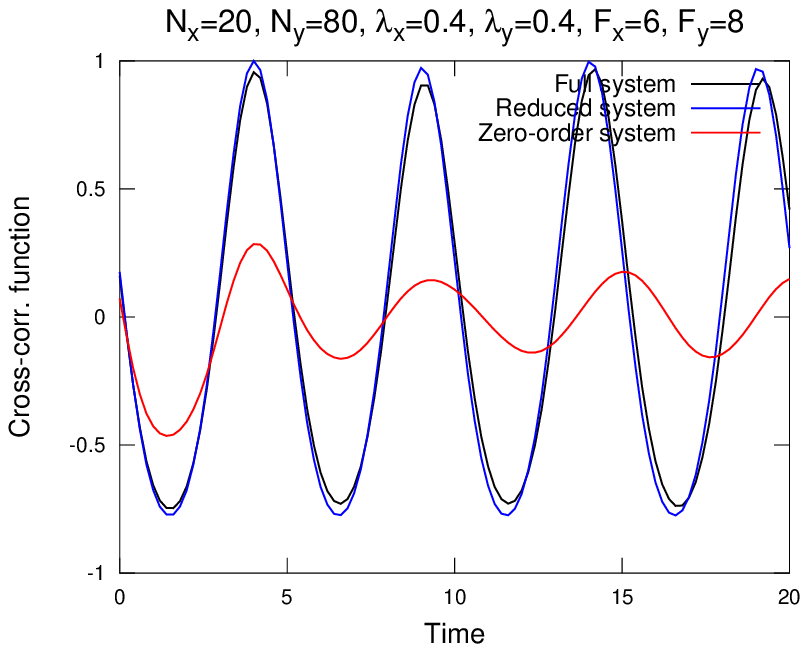}%
\picturehere{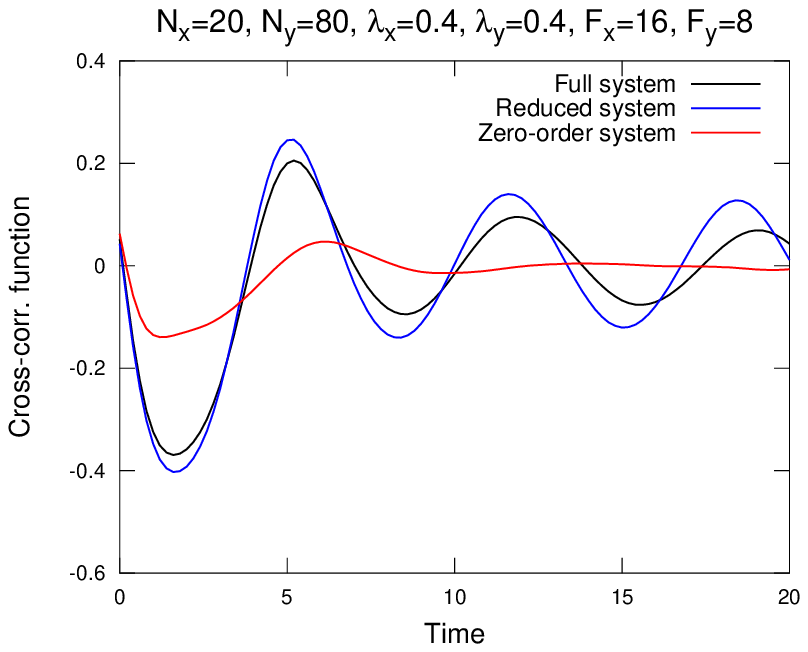}\\%
\picturehere{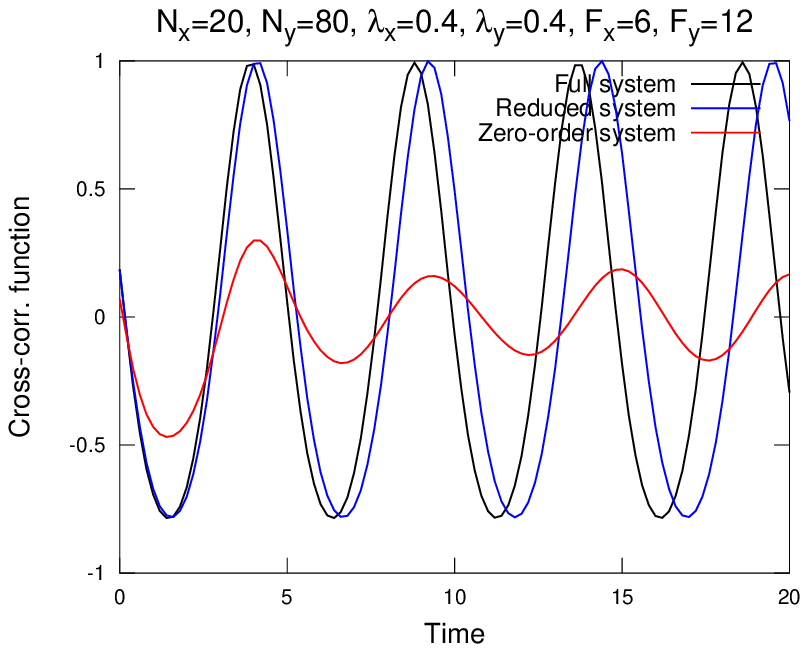}%
\picturehere{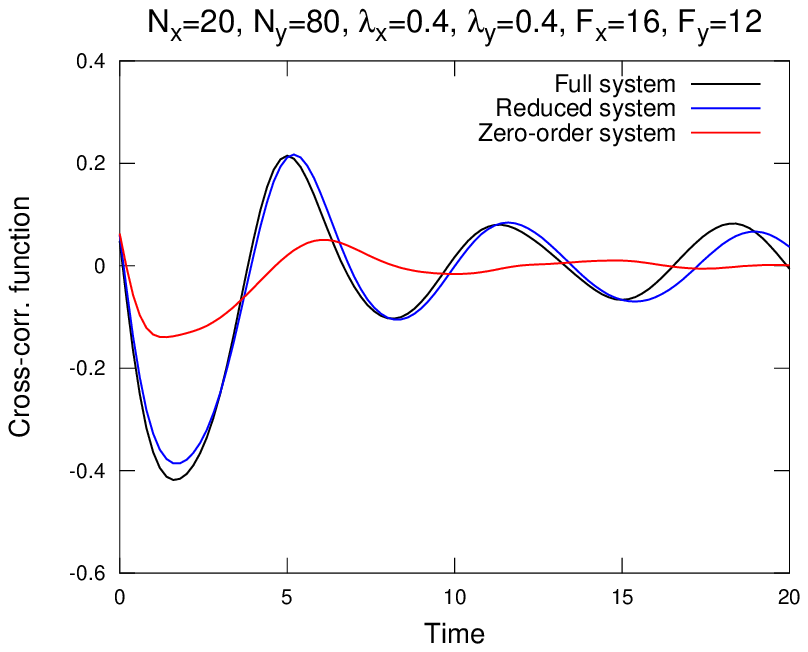}%
\caption{Time cross-correlation functions, $\lambda_x=\lambda_y=0.4$.}
\label{fig:ccorr_l04}
\end{figure}

In Figures \ref{fig:ccorr_l03} and \ref{fig:ccorr_l04} we show the
time cross-correlation functions of the slow dynamics for the full
two-scale Lorenz model, the reduced closed model for the slow
variables alone in \eqref{eq:dyn_sys_slow_limiting_x_lor}, and its
poor man's zero order version without the linear correction
term. Observe that for the more weakly coupled regimes with
$\lambda_x=\lambda_y=0.3$ the time cross-correlation functions look
similar, however, it is seen that the reduced model with the
correction term reproduces the time cross-correlation functions more
precisely than the zero-order model. In the more strongly coupled
regime with $\lambda_x=\lambda_y=0.4$, the reduced model in
\eqref{eq:dyn_sys_slow_limiting_x_lor} becomes much more precise than
its poor man's zero-order version: here, for the weakly chaotic regime
with $F_x=6$ the time cross-correlation functions of the full
two-scale Lorenz model do not exhibit decay (indicating very weak
mixing), and the reduced model in
\eqref{eq:dyn_sys_slow_limiting_x_lor} reproduces the
cross-correlation functions of the full two-scale Lorenz model rather
well, while its zero-order version fails.
In addition, in Table \ref{tab:ccorrs} we show the $L_2$-errors in
time cross-correlation functions (for the correlation time interval of
20 time units, as in Figures \ref{fig:ccorr_l03} and
\ref{fig:ccorr_l04}) between the full two-scale Lorenz model and the
two reduced models. Observe that, generally, the reduced system with
the linear correction term in \eqref{eq:dyn_sys_slow_limiting_x_lor}
produces more precise results than its poor man's version without the
correction term.
\begin{table}
\begin{center}
\begin{tabular}{|c|}
\hline
\begin{tabular}{c|c}
$\lambda_{x,y}=0.3$, $F_y=8$ & $\lambda_{x,y}=0.3$, $F_y=12$ \\
\hline
\begin{tabular}{c||c|c}
 & Av. & Z.O. \\
\hline
$F_x=6$ & $6.825\cdot 10^{-2}$ & $0.1437$ \\
$F_x=16$ & $0.1094$ & $0.2134$ \\
\end{tabular}
&
\begin{tabular}{c||c|c}
 & Av. & Z.O. \\
\hline
$F_x=6$ & $6.838\cdot 10^{-2}$ & $0.1799$ \\
$F_x=16$ & $3.687\cdot 10^{-2}$ & $0.2548$ \\
\end{tabular}\\
\end{tabular}\\
\hline
\begin{tabular}{c|c}
$\lambda_{x,y}=0.4$, $F_y=8$ & $\lambda_{x,y}=0.4$, $F_y=12$ \\
\hline
\begin{tabular}{c||c|c}
 & Av. & Z.O. \\
\hline
$F_x=6$ & $5.313\cdot 10^{-2}$ & $0.3666$ \\
$F_x=16$ & $0.1258$ & $0.3232$ \\
\end{tabular}
&
\begin{tabular}{c||c|c}
 & Av. & Z.O. \\
\hline
$F_x=6$ & $0.3137$ & $0.3942$ \\
$F_x=16$ & $6.953\cdot 10^{-2}$ & $0.3321$ \\
\end{tabular}\\
\end{tabular}\\
\hline
\end{tabular}
\end{center}
\caption{$L_2$-errors between the time cross-correlation functions of
  the slow variables of the full two-scale Lorenz model and the two
  reduced models. Notations: ``Red.'' stands for ``Reduced'' (that is,
  \eqref{eq:dyn_sys_slow_limiting_x_lor}), and ``Z.O.'' stands for
  ``Zero-order'', the poor man's version of
  \eqref{eq:dyn_sys_slow_limiting_x_lor}.}
\label{tab:ccorrs}
\end{table}

\subsection{Energy autocorrelation functions}

\begin{figure}
\picturehere{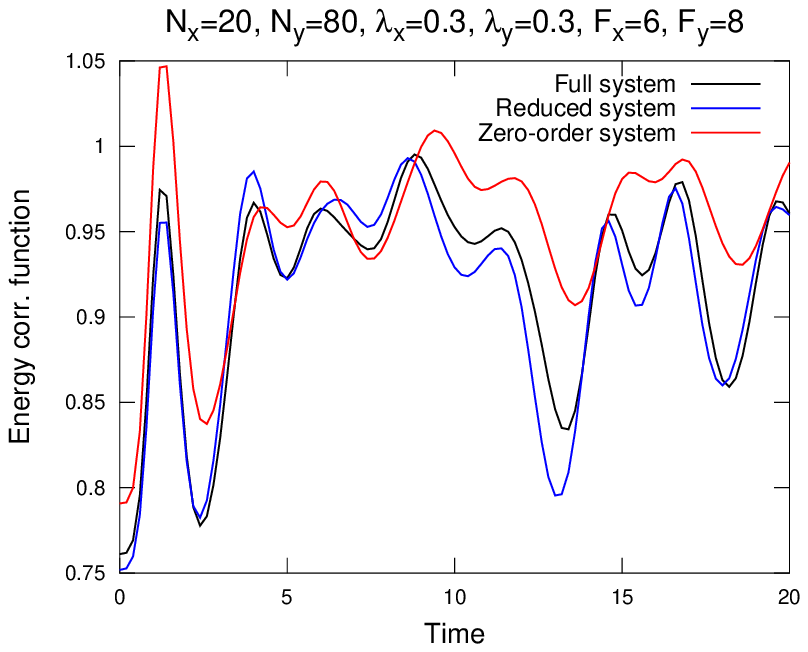}%
\picturehere{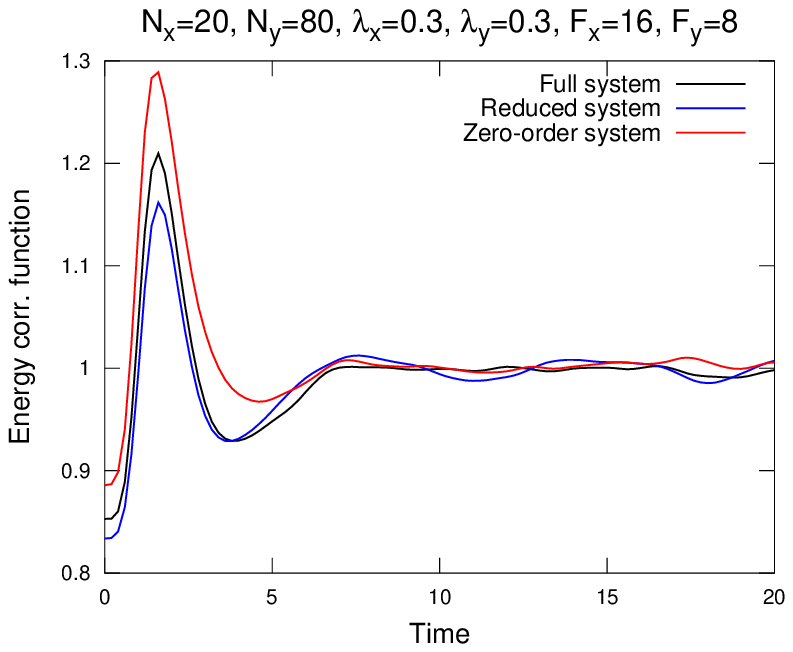}\\%
\picturehere{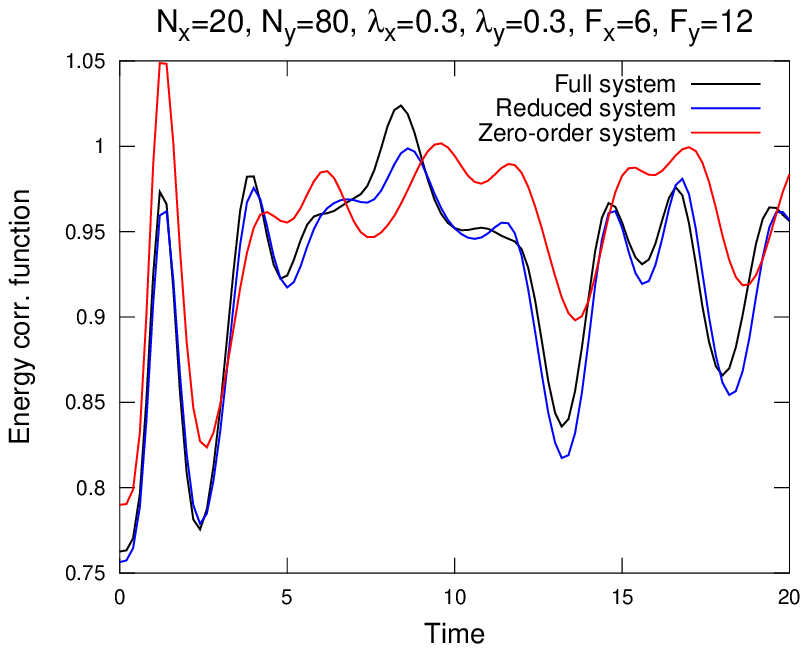}%
\picturehere{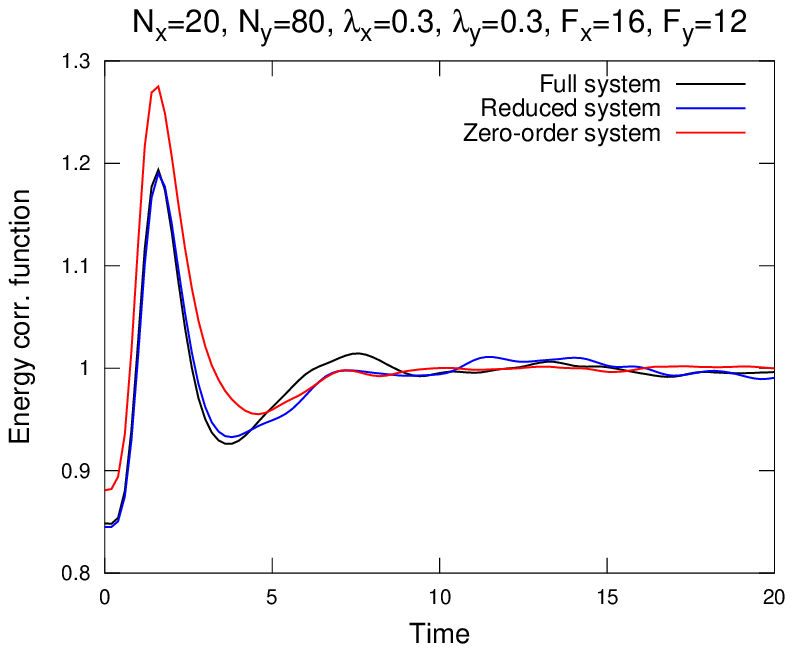}%
\caption{Energy autocorrelation functions, $\lambda_x=\lambda_y=0.3$.}
\label{fig:Kcorr_l03}
\end{figure}

\begin{figure}
\picturehere{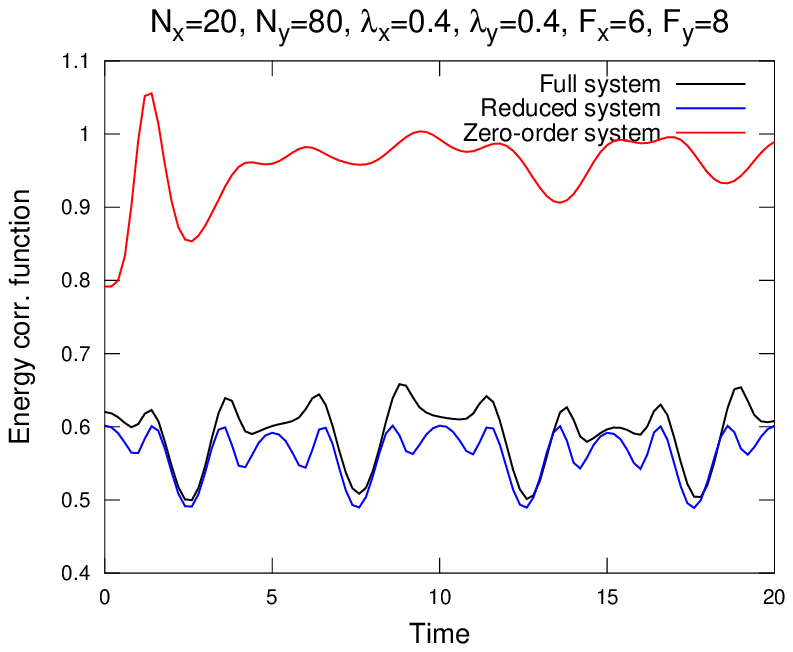}%
\picturehere{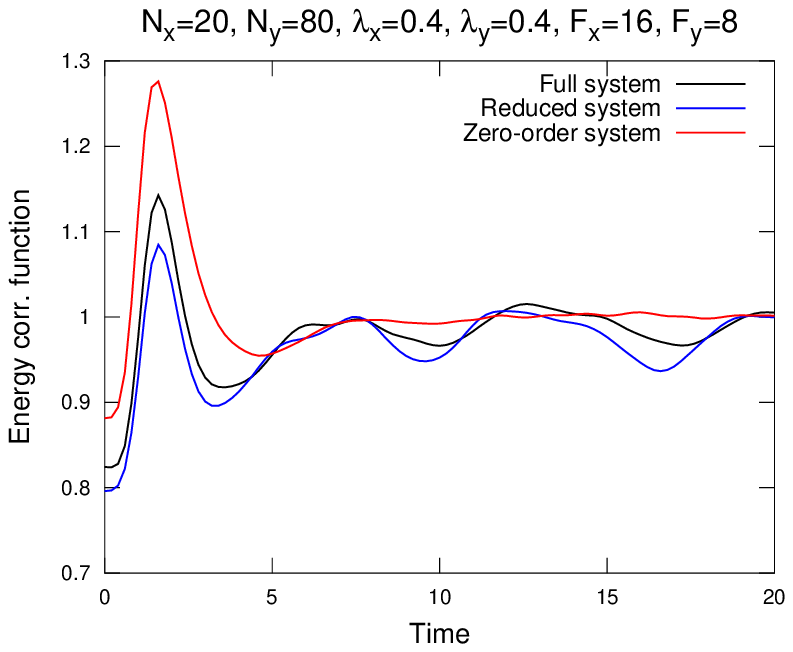}\\%
\picturehere{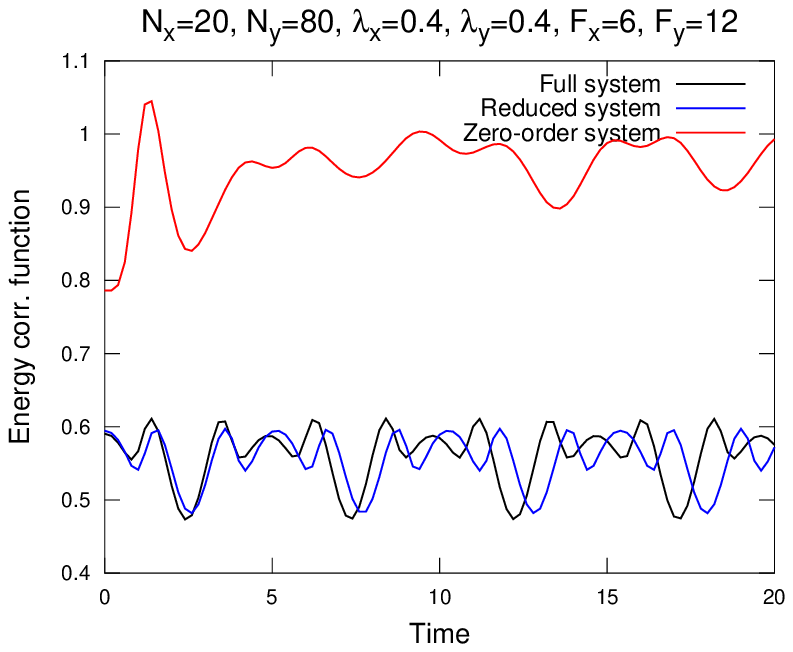}%
\picturehere{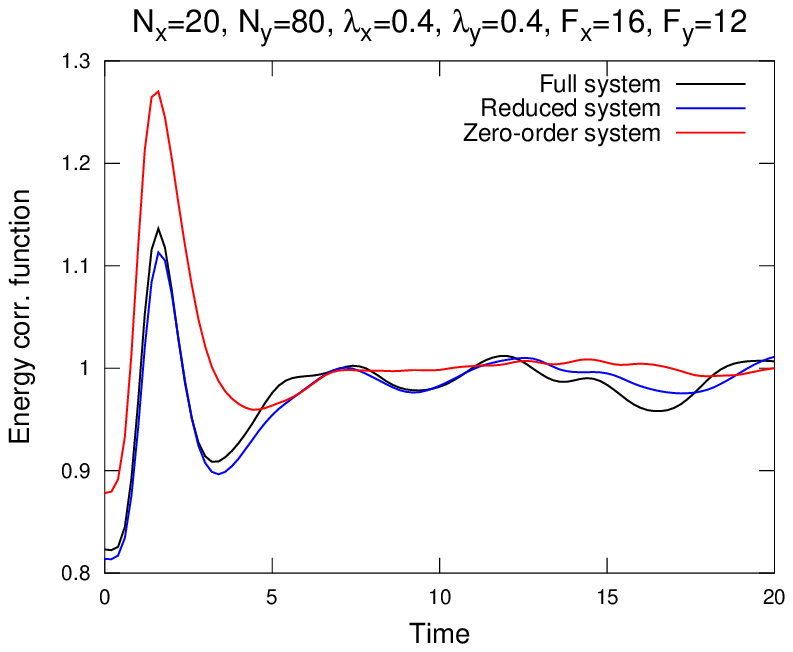}%
\caption{Energy autocorrelation functions, $\lambda_x=\lambda_y=0.4$.}
\label{fig:Kcorr_l04}
\end{figure}

In Figures \ref{fig:Kcorr_l03} and \ref{fig:Kcorr_l04} we show the
energy autocorrelation functions of the slow dynamics for the full
two-scale Lorenz model, the reduced closed model for the slow
variables alone in \eqref{eq:dyn_sys_slow_limiting_x_lor}, and its
poor man's zero order version without the linear correction
term. Observe that for the more weakly coupled regimes with
$\lambda_x=\lambda_y=0.3$ the energy autocorrelation functions look
similar, although the reduced model with the correction term
reproduces the energy autocorrelation functions more precisely than
the zero-order model. In the more strongly coupled regime with
$\lambda_x=\lambda_y=0.4$, the reduced model in
\eqref{eq:dyn_sys_slow_limiting_x_lor} is much more precise than its
poor man's zero-order version: here, for the weakly chaotic regime
with $F_x=6$ the energy autocorrelation functions of the full
two-scale Lorenz model is significantly sub-Gaussian, and the reduced
model in \eqref{eq:dyn_sys_slow_limiting_x_lor} reproduces the
sub-Gaussianity of the energy autocorrelation functions of the full
two-scale Lorenz model rather well, while its zero-order version
fails.
In addition, in Table \ref{tab:Kcorrs} we show the $L_2$-errors in
energy autocorrelation functions (for the correlation time interval of
20 time units, as in Figures \ref{fig:Kcorr_l03} and
\ref{fig:Kcorr_l04}) between the full two-scale Lorenz model and the
two reduced models. Observe that, generally, the reduced system with
the linear correction term in \eqref{eq:dyn_sys_slow_limiting_x_lor}
produces more precise results than its poor man's version without the
correction term.
\begin{table}
\begin{center}
\begin{tabular}{|c|}
\hline
\begin{tabular}{c|c}
$\lambda_{x,y}=0.3$, $F_y=8$ & $\lambda_{x,y}=0.3$, $F_y=12$ \\
\hline
\begin{tabular}{c||c|c}
 & Av. & Z.O. \\
\hline
$F_x=6$ & $8.911\cdot 10^{-3}$ & $2.027\cdot 10^{-2}$ \\
$F_x=16$ & $6.885\cdot 10^{-3}$ & $1.434\cdot 10^{-2}$ \\
\end{tabular}
&
\begin{tabular}{c||c|c}
 & Av. & Z.O. \\
\hline
$F_x=6$ & $6.783\cdot 10^{-3}$ & $2.131\cdot 10^{-2}$ \\
$F_x=16$ & $4.154\cdot 10^{-3}$ & $1.455\cdot 10^{-2}$ \\
\end{tabular}\\
\end{tabular}\\
\hline
\begin{tabular}{c|c}
$\lambda_{x,y}=0.4$, $F_y=8$ & $\lambda_{x,y}=0.4$, $F_y=12$ \\
\hline
\begin{tabular}{c||c|c}
 & Av. & Z.O. \\
\hline
$F_x=6$ & $2.66\cdot 10^{-2}$ & $0.2779$ \\
$F_x=16$ & $9.746\cdot 10^{-3}$ & $2.284\cdot 10^{-2}$ \\
\end{tabular}
&
\begin{tabular}{c||c|c}
 & Av. & Z.O. \\
\hline
$F_x=6$ & $3.499\cdot 10^{-2}$ & $0.3125$ \\
$F_x=16$ & $5.49\cdot 10^{-3}$ & $2.414\cdot 10^{-2}$ \\
\end{tabular}\\
\end{tabular}\\
\hline
\end{tabular}
\end{center}
\caption{$L_2$-errors between the energy autocorrelation functions of
  the slow variables of the full two-scale Lorenz model and the two
  reduced models. Notations: ``Red.'' stands for ``Reduced'' (that is,
  \eqref{eq:dyn_sys_slow_limiting_x_lor}), and ``Z.O.'' stands for
  ``Zero-order'', the poor man's version of
  \eqref{eq:dyn_sys_slow_limiting_x_lor}.}
\label{tab:Kcorrs}
\end{table}

\section{Summary}
\label{sec:sum}

In this work we develop a simple method of constructing the closed
reduced model for slow variables of a multiscale model with linear
coupling, which requires only a single computation of the mean state
and the time autocorrelation function for the fast dynamics with a
fixed state of the slow variables, which is located in the region
where the slow dynamics evolve (here, the mean state of the slow
dynamics is used). The method is based on the first-order Taylor
expansion of the averaged coupling term with respect to the slow
variables, which is computed using the linear fluctuation-dissipation
theorem. We demonstrate through the computations with the
appropriately rescaled two-scale Lorenz 96 model \cite{Abr8} that,
with simple linear coupling in both slow and fast variables, the
developed reduced model produces quite comparable statistics to what
is exhibited by the complete two-scale Lorenz model. Below we outline
the main advantages of the new method:
\begin{itemize}
\item The reduced model is simple. It requires only the mean state,
  covariance and the time autocorrelation function for the fast
  variables, computed for a single fixed state of the slow
  variables. Since only the statistics of the time series of the fast
  dynamics are needed, the structure of the right-hand side of the
  equations for the fast variables need not be known explicitly --
  this part of dynamics can be provided as a ``black box''. Also, the
  structure of the nonlinear $\BS x$-dependent part of the right-hand
  side of the equations for the slow variables need not be known to
  construct the approximation, and existing computational routines can
  be used for it. The only correction in the forward time-stepping
  routine is the linear correction term.
\item The reduced model is {\em a priori}. It lacks parameters which
  have to be adjusted {\em a posteriori} to ``fit'' the statistical
  properties of the full multiscale dynamics. In fact, statistical
  properties of the full multiscale dynamics need not be known to
  construct the reduced model (although certain statistics of the fast
  variables with an appropriate fixed slow state need to be computed).
\item The reduced model is parsimonious. It requires only a simple
  linear correction to achieve consistently better performance than
  that of a corresponding zero-order model.
\item The reduced model is practical. It can be implemented even when
  the statistics for the slow variables of the complete multiscale
  dynamics cannot be obtained due to its computational complexity
  (although a rough estimate of the mean state, or a nearby state is
  needed), which makes the approach potentially suitable for
  comprehensive global circulation models in geophysics. Additionally,
  existing zero-order models (such as the T21 barotropic model
  \cite{AbrMaj6,Fra,Sel}) can be retrofitted with the linear
  correction term.
\end{itemize}
In the future work, the author intends to collaborate with
geophysicists to create more realistic reduced models for geophysical
dynamics, including retrofitting existing closed models for slow
dynamics with the linear correction term.

\begin{acknowledgment}
The author thanks Ibrahim Fatkullin, Ilya Timofeyev and Gregor Kova\v
ci\v c for fruitful discussions. The author is supported by the
National Science Foundation CAREER grant DMS-0845760, and the Office
of Naval Research grants N00014-09-0083 and 25-74200-F6607.
\end{acknowledgment}


\begin{thebibliography}{10}

\bibitem{Abr5}
R.~Abramov.
\newblock Short-time linear response with reduced-rank tangent map.
\newblock {\em Chin. Ann. Math.}, 30B(5):447--462, 2009.

\bibitem{Abr6}
R.~Abramov.
\newblock Approximate linear response for slow variables of deterministic or
  stochastic dynamics with time scale separation.
\newblock {\em J. Comput. Phys.}, 229(20):7739--7746, 2010.

\bibitem{Abr7}
R.~Abramov.
\newblock Improved linear response for stochastically driven systems.
\newblock {\em Front. Math. China}, 2011.
\newblock submitted.

\bibitem{Abr8}
R.~Abramov.
\newblock Suppression of chaos at slow variables by rapidly mixing fast
  dynamics through linear energy-preserving coupling.
\newblock {\em Comm. Math. Sci.}, 2011.
\newblock submitted.

\bibitem{AbrMaj5}
R.~Abramov and A.J. Majda.
\newblock Blended response algorithms for linear fluctuation-dissipation for
  complex nonlinear dynamical systems.
\newblock {\em Nonlinearity}, 20:2793--2821, 2007.

\bibitem{AbrMaj4}
R.~Abramov and A.J. Majda.
\newblock New approximations and tests of linear fluctuation-response for
  chaotic nonlinear forced-dissipative dynamical systems.
\newblock {\em J. Nonlin. Sci.}, 18(3):303--341, 2008.

\bibitem{AbrMaj6}
R.~Abramov and A.J. Majda.
\newblock New algorithms for low frequency climate response.
\newblock {\em J. Atmos. Sci.}, 66:286--309, 2009.

\bibitem{AbrMaj7}
R.~Abramov and A.J. Majda.
\newblock Low frequency climate response of quasigeostrophic wind-driven ocean
  circulation.
\newblock {\em J. Phys. Oceanogr.}, 2011.
\newblock submitted.

\bibitem{BuiMilPal}
R.~Buizza, M.~Miller, and T.~Palmer.
\newblock Stochastic representation of model uncertainty in the {ECMWF}
  {E}nsemble {P}rediction {S}ystem.
\newblock {\em Q. J. R. Meteor. Soc.}, 125:2887--2908, 1999.

\bibitem{CroVan}
D.~Crommelin and E.~Vanden-Eijnden.
\newblock Subgrid scale parameterization with conditional {M}arkov chains.
\newblock {\em J. Atmos. Sci.}, 65:2661--2675, 2008.

\bibitem{ELi}
W.~E and X.~Li.
\newblock Some recent progress in multiscale modeling.
\newblock In S.~Attinger and P.~Koumoutsakos, editors, {\em Multiscale
  Modelling and Simulation}, volume~39 of {\em LNCSE}. Springer, 2004.

\bibitem{EckRue}
J.~Eckmann and D.~Ruelle.
\newblock Ergodic theory of chaos and strange attractors.
\newblock {\em Rev. Mod. Phys.}, 57(3):617--656, 1985.

\bibitem{FatVan}
I.~Fatkullin and E.~Vanden-Eijnden.
\newblock A computational strategy for multiscale systems with applications to
  {L}orenz 96 model.
\newblock {\em J. Comp. Phys.}, 200:605--638, 2004.

\bibitem{Fra}
C.~Franzke.
\newblock Dynamics of low-frequency variability: {B}arotropic mode.
\newblock {\em J. Atmos. Sci.}, 59:2909--2897, 2002.

\bibitem{FraMajVan}
C.~Franzke, A.J. Majda, and E.~Vanden-Eijnden.
\newblock Low-order stochastic model reduction for a realistic barotropic model
  climate.
\newblock {\em J. Atmos. Sci.}, 62:1722--1745, 2005.

\bibitem{Has}
K.~Hasselmann.
\newblock Stochastic climate models, part {I}, theory.
\newblock {\em Tellus}, 28:473--485, 1976.

\bibitem{Lor}
E.~Lorenz.
\newblock Predictability: A problem partly solved.
\newblock In {\em Proceedings of the Seminar on Predictability}, Shinfield
  Park, Reading, England, 1996. ECMWF.

\bibitem{LorEma}
E.~Lorenz and K.~Emanuel.
\newblock Optimal sites for supplementary weather observations.
\newblock {\em J. Atmos. Sci.}, 55:399--414, 1998.

\bibitem{MajAbrGro}
A.J. Majda, R.~Abramov, and M.~Grote.
\newblock {\em Information Theory and Stochastics for Multiscale Nonlinear
  Systems}, volume~25 of {\em CRM Monograph Series of Centre de Recherches
  Math\'ematiques, Universit\'e de Montr\'eal}.
\newblock American Mathematical Society, 2005.
\newblock {ISBN} 0-8218-3843-1.

\bibitem{MajTimVan}
A.J. Majda, I.~Timofeyev, and E.~Vanden-Eijnden.
\newblock Models for stochastic climate prediction.
\newblock {\em Proc. Natl. Acad. Sci.}, 96:14687--14691, 1999.

\bibitem{MajTimVan3}
A.J. Majda, I.~Timofeyev, and E.~Vanden-Eijnden.
\newblock A mathematical framework for stochastic climate models.
\newblock {\em Comm. Pure Appl. Math.}, 54:891--974, 2001.

\bibitem{MajTimVan4}
A.J. Majda, I.~Timofeyev, and E.~Vanden-Eijnden.
\newblock A priori tests of a stochastic mode reduction strategy.
\newblock {\em Physica D}, 170:206--252, 2002.

\bibitem{MajTimVan2}
A.J. Majda, I.~Timofeyev, and E.~Vanden-Eijnden.
\newblock Systematic strategies for stochastic mode reduction in climate.
\newblock {\em J. Atmos. Sci.}, 60:1705--1722, 2003.

\bibitem{Pal3}
T.~Palmer.
\newblock A nonlinear dynamical perspective on model error: A proposal for
  nonlocal stochastic-dynamic parameterization in weather and climate
  prediction models.
\newblock {\em Q. J. R. Meteor. Soc.}, 127:279--304, 2001.

\bibitem{Pap}
G.~Papanicolaou.
\newblock Introduction to the asymptotic analysis of stochastic equations.
\newblock In R.~DiPrima, editor, {\em Modern modeling of continuum phenomena},
  volume~16 of {\em Lectures in Applied Mathematics}. American Mathematical
  Society, 1977.

\bibitem{Ris}
F.~Risken.
\newblock {\em The {F}okker-{P}lanck Equation}.
\newblock Springer-Verlag, New York, second edition, 1988.

\bibitem{Rue3}
D.~Ruelle.
\newblock A measure associated with {A}xiom {A} attractors.
\newblock {\em Amer. J. Math.}, 98:619--654, 1976.

\bibitem{Rue1}
D.~Ruelle.
\newblock Differentiation of {SRB} states.
\newblock {\em Comm. Math. Phys.}, 187:227--241, 1997.

\bibitem{Rue2}
D.~Ruelle.
\newblock General linear response formula in statistical mechanics, and the
  fluctuation-dissipation theorem far from equilibrium.
\newblock {\em Phys. Lett. A}, 245:220--224, 1998.

\bibitem{Sel}
F.~Selten.
\newblock An efficient description of the dynamics of barotropic flow.
\newblock {\em J. Atmos. Sci.}, 52:915--936, 1995.

\bibitem{OrnUhl}
G.~Uhlenbeck and L.~Ornstein.
\newblock On the theory of the {B}rownian motion.
\newblock {\em Phys. Rev.}, 36:823--841, 1930.

\bibitem{Van}
E.~Vanden-Eijnden.
\newblock Numerical techniques for multiscale dynamical systems with stochastic
  effects.
\newblock {\em Comm. Math. Sci.}, 1:385--391, 2003.

\bibitem{Vol}
V.~Volosov.
\newblock Averaging in systems of ordinary differential equations.
\newblock {\em Russian Math. Surveys}, 17:1--126, 1962.

\bibitem{Wilks}
D.~Wilks.
\newblock Effects of stochastic parameterizations in the {L}orenz '96 system.
\newblock {\em Q. J. R. Meteorol. Soc.}, 131:389--407, 2005.

\bibitem{You}
L.-S. Young.
\newblock What are {SRB} measures, and which dynamical systems have them?
\newblock {\em J. Stat. Phys.}, 108(5-6):733--754, 2002.

\end{thebibliography}
\end{document}